\documentclass[12pt,reqno,twoside,fleqn]{amsart}
\usepackage{amscd,amsmath,amssymb,amsthm,palatino}%
\usepackage[mathcal]{euler}

\numberwithin{equation}{section}
\makeatletter
\renewcommand{\section}{\@startsection {section}{1}{\z@}%
                                   {-3.5ex \@plus -1ex \@minus -.2ex}%
                                   {.5\linespacing}%
                                   {\normalfont\scshape\centering}}
\makeatother

\newtheorem{thm}{Theorem}[section]
\newtheorem{lem}[thm]{Lemma}
\newtheorem{cor}[thm]{Corollary}
\newtheorem{prop}[thm]{Proposition}
  
\theoremstyle{definition}
\newtheorem{definition}{Definition}[section]

\theoremstyle{remark}
\newtheorem*{example}{Example}

\def\beq#1\eeq{\begin{equation}#1\end{equation}}
\makeatletter
\@ifpackageloaded{euler}{
 
 \newcommand{\onto}{\to\mkern-14mu\to}
 \def\rightarrowfill@#1{\m@th\setboxz@h{$#1\relbar$}\ht\z@\z@
   $#1\copy\z@\mkern-6mu\cleaders
   \hbox{$#1\mkern-2mu\box\z@\mkern-2mu$}\hfill
   \mkern-6mu\mathord\rightarrow$}
 \def\leftarrowfill@#1{\m@th\setboxz@h{$#1\relbar$}\ht\z@\z@
   $#1\mathord\leftarrow\mkern-6mu\cleaders
   \hbox{$#1\mkern-2mu\copy\z@\mkern-2mu$}\hfill
   \mkern-6mu\box\z@$}
 \def\B@R#1#2{\raisebox{-.07ex}{$#1#2$}\mkern-6mu}
 \renewcommand{\hbar}{{\mspace{1mu}\mathpalette\B@R{\mathchar'26}h}}
 \newcommand{\Wedge}{\mathord{\wedge}}
 \newcommand{\espace}{\hspace{.1pt}}
}{
 
 \DeclareMathSymbol{\onto}{\mathrel}{AMSa}{"10}
 \renewcommand{\hbar}{{\mathchar'26\mkern-9muh}}
 \newcommand{\Wedge}{\mathsf{\Lambda}}
 \newcommand{\espace}{}
}
\makeatother
\newcommand{\C}{\mathcal{C}}
\newcommand{\co}{\mathbb{C}}
\newcommand{\cs}{\mbox{\upshape C}\ensuremath{{}^*}} 
\newcommand{\csr}{{\mathrm C}^*_{\mathrm r}}

\newcommand{\R}{\mathbb{R}}
\newcommand{\Z}{\mathbb{Z}} 
\newcommand{\into}{\hookrightarrow}
\newcommand{\isom}{\mathrel{\widetilde\longrightarrow}}

\DeclareMathOperator{\rk}{rk}
\newcommand{\inner}{\mathbin{\raise1.5pt\hbox{$\lrcorner$}}}
\newcommand{\binner}{\mathbin{\raise1.5pt\hbox{$\llcorner$}}}

\newcommand{\TcM}{T_{\co}M}
\newcommand{\K}{\mathcal{K}}
\newcommand{\G}{\mathcal{G}}
\DeclareMathOperator{\pr}{pr}
\newcommand{\T}{\mathbb{T}}
\DeclareMathOperator{\inv}{inv}
\newcommand{\Mt}{\widetilde M}
\newcommand{\E}{\mathcal{E}}
\newcommand{\D}{\mathcal{D}}
\newcommand{\F}{\mathcal{F}}
\renewcommand{\P}{\mathcal{P}}
\newcommand{\Q}{\mathcal{Q}}
\newcommand{\Lie}{\mathcal{L}}
\renewcommand{\Im}{\mathop{\mathrm{Im}}}
\newcommand{\s}{\mathsf{s}}
\renewcommand{\t}{\mathsf{t}}
\newcommand{\m}{\mathsf{m}}
\newcommand{\unit}{\mathsf{1}}

\newcommand{\abs}[1]{\lvert#1\rvert}

\newcommand{\Kahler}{K\"ahler} 
\newcommand{\Wedgem}{\Wedge^{\mathrm{max}}}
\DeclareMathOperator{\Diff}{Diff}
\newcommand{\tf}{\mathfrak{t}}
\newcommand{\X}{\mathcal{X}}
\newcommand{\h}{\mathfrak{H}}

\newcommand{\Alg}{\mathsf{A}}
\newcommand{\p}{\mathfrak{P}}
\newcommand{\Grp}{\mathsf{G}}
\newcommand{\Sig}{{\sf\Sigma}}
\newcommand{\PI}{{\sf\Pi}}
\DeclareMathOperator{\Pair}{Pair}
\DeclareMathOperator{\Tw}{Tw}
\newcommand{\BS}{\text{B-S}}
\newcommand{\A}{\mathcal{A}}
\DeclareMathOperator{\curv}{curv}

\newcommand{\action}{\rho}
\newcommand{\Action}{R}
\newcommand{\daction}{\tau}
\newcommand{\auxiliary}{B}

\title{A Groupoid Approach to Quantization}
\author{Eli Hawkins}

\begin{document}
\maketitle
\begin{center}
\vspace{-4ex}
\emph{\small Institute for Mathematics, Astrophysics, and Particle Physics}\\
\emph{\small Radboud University Nijmegen, The Netherlands}\\
{\small mrmuon@mac.com}\\
\end{center}

\begin{abstract}
Many interesting \cs-algebras can be viewed as quantizations of Poisson manifolds. I propose that a Poisson manifold may be quantized by a twisted polarized convolution \cs-algebra of a symplectic groupoid. Toward this end, I define polarizations for Lie groupoids and sketch the construction of this algebra. A large number of examples show that this idea unifies previous geometric constructions, including geometric quantization of symplectic manifolds and the \cs-algebra of a Lie groupoid.
\end{abstract}

\section{Introduction}
Many interesting \cs-algebras can be regarded as quantizations. These include the algebra of compact operators on a Hilbert space, the \cs-algebra of a Lie group \cite{rie4}, the noncommutative torus \cite{rie2}, the crossed product \cs-algebra of a group acting on a manifold, the \cs-algebra of a foliation, the \cs-algebra of a Lie groupoid \cite{l-r}, quantum groups, and of course any of the algebras appearing in quantum physics. 
A better understanding of quantization should provide many more examples and lead to new tools for understanding these algebras.

The space of continuous functions on a manifold is a commutative algebra under the operations of pointwise addition and multiplication. A bivector field on a smooth manifold determines an antisymmetric bracket of differentiable functions. If this satisfies the Jacobi identity, then it is a \emph{Poisson bracket} and can be treated as a first order correction to multiplication. A \emph{Poisson manifold} \cite{vai2} is a manifold with such a bivector field.

In this way, a Poisson manifold may be regarded as a geometrical approximation to a noncommutative algebra. I am specifically interested in \cs-algebras here. I will say that a \cs-algebra $\A$ \emph{quantizes} a Poisson manifold $M$ if the Poisson algebra of functions on $M$ approximates $\A$. This is in the general \cs-algebraic approach to quantization advocated by Rieffel, but I will clarify the definition below. This paper is concerned with the problem of constructing a quantization systematically from a Poisson manifold with some additional structure. 

\subsection{Symplectic Groupoids and the Dictionary}
Symplectic groupoids were invented independently by Karas\"ev, Weinstein, and Zakrzewski as a tool for quantization. The base manifold of a symplectic groupoid is a Poisson manifold, and when it exists, the symplectic groupoid of a Poisson manifold is unique modulo covering and isomorphism. The Poisson structure can be thought of as an infinitesimal structure that is ``integrated'' by the groupoid. A symplectic groupoid is supposed to be an intermediate structure on the way to quantizing its base Poisson manifold.

A specific strategy for quantizing in this way was outlined by Weinstein \cite{wei3,wei4,wei5}. His approach is to mainly regard the groupoid as a symplectic manifold with some additional structures and to apply geometric quantization in a particular sense known as \emph{the dictionary}. This is a proposed correspondence between geometric (classical) and algebraic (quantum) objects and constructions:
\begin{center}
\begin{tabular}{|c|c|}
\hline
\emph{Geometrical} & \emph{Algebraic} \\
\hline
Symplectic manifold $M$ & Vector space \\
$M^-$ (with opposite symplectic form) & Dual vector space \\
Cartesian product & Tensor product \\
Lagrangian submanifold & Vector\\
\hline
\end{tabular}
\end{center}

Let $\Sigma$ be a symplectic groupoid; the corresponding vector space is supposed to become the algebra $\A$. The graph of multiplication is a Lagrangian submanifold of $\Sigma^-\times\Sigma^-\times\Sigma$; this should correspond to a vector in $\A^*\otimes\A^*\otimes\A$ --- or equivalently, a bilinear map from $\A\times\A$ to $\A$; this is supposed to be the product. The unit manifold of $\Sigma$ is a Lagrangian submanifold; this should correspond to an element of $\A$ which is supposed to be the unit.

This approach has only been successfully carried out for a few simple examples; see \cite{g-v,wei3,wei4}.
The problem is that the product is only guaranteed to be associative if the dictionary can be implemented exactly. In practice, most things in geometric quantization are only approximately true. They become exact in the classical limit.

Part of the problem is that the vector space cannot be constructed from a symplectic manifold alone. It requires additional structure including a prequantization and polarization. The dictionary-based approach does not assume the polarization to be compatible with the groupoid structures in any way.

\subsection{Alternative Approach}
The purpose of this paper is to propose a different approach to quantizing Poisson manifolds using symplectic groupoids. Rather than regarding $\Sigma$ primarily as a symplectic manifold, I treat it primarily as a groupoid. I propose that quantization is a modification of the construction of the convolution \cs-algebra of the groupoid $\Sigma$.

With this in mind, I propose a definition (Def.~\ref{polarization}) of polarization of a groupoid in parallel with the notion of polarization of a symplectic manifold. The prototypical groupoid polarization is a simple foliation of a groupoid such that the leaf space is also a groupoid. In some ways, I will treat an arbitrary polarized groupoid as a stand-in for a quotient groupoid that may not exist.

My proposed recipe for geometric quantization of a Poisson manifold $M$ consists of the following steps: 
\begin{enumerate}
\item
Construct an $\s$-connected symplectic groupoid over $M$.
\item
Construct a prequantization $(\sigma,L,\nabla)$ of $\Sigma$.
\item
Choose a symplectic groupoid polarization $\P$ of $\Sigma$. That is, some $\P$ satisfying both the definitions of a symplectic polarization and a groupoid polarization.
\item
Construct a ``half-form'' bundle (or sheaf) $\Omega^{1/2}_\P$.
\item
$M$ is quantized by the twisted, polarized convolution algebra $\cs_\P(\Sigma,\sigma)$. This is essentially a convolution algebra of sections of $L$ that are polarized by $\P$.
\end{enumerate}
The first four steps each entail existence and uniqueness issues. The required structure may not exist and is not necessarily unique.

My definition of a symplectic groupoid  polarization is very restrictive. Such polarizations may not exist in sufficient generality to quantize all Poisson manifolds that should be quantizable.
Nevertheless, this is the optimal scenario. If the symplectic groupoid approach to quantization doesn't work with this type of polarization, then it will surely not work for more general polarizations.

I will show in examples that such polarizations actually do exist in many cases. Not only that, but these examples appear to reproduce \emph{every} example of geometrically constructed quantization that I am aware of.

\subsection{Quantization}
\label{Quantization}
The term ``quantization'' has \emph{many} meanings. Planck's original ``quantum hypothesis'' postulated discrete units of energy. As the understanding of this physics developed, it proved to be more about noncommutativity than discreteness. Canonical quantization emerged as the process of constructing a quantum mechanical model from a classical one, using the canonical commutation relations between position and momentum.

Unfortunately, the canonical commutation relations are not coordinate independent, and canonical quantization suffers from operator ordering ambiguities. For the most part, these problems are not seriously troubling to physicists, but they have led to interesting mathematical abstractions, including geometric quantization \cite{a-e,sni,woo} and formal deformation quantization \cite{bffls}. 

Although the phase space of a classical system is a symplectic manifold, formal deformation quantization applies to any Poisson manifold. A formal deformation quantization of a Poisson manifold is an associative product on the space $\C^\infty(M)[[\hbar]]$ of smooth functions valued in formal power series in ``Planck's constant'' $\hbar$. The product is required to be pointwise multiplication modulo $\hbar$, and the commutator is required to be $i\hbar$ times the Poisson bracket, modulo $\hbar^2$.

Because formal deformation quantization uses formal power series, it is not possible to insert a specific value for $\hbar$. This is an unfortunate departure from physics, where $\hbar$ is Planck's \emph{constant} and conventionally set equal to $1$ by many theoreticians.

A more concrete approach to quantization was initiated by Rieffel \cite{rie6}. This was originally stated in terms of an $\hbar$-dependent product, but it can be restated in terms of a continuous field of \cs-algebras. A \emph{strict deformation quantization} of a Poisson manifold $M$ consists of a continuous field of \cs-algebras $\{\A_\hbar\}$ and quantization maps $Q_\hbar : \C^\infty_0(M)\to \A_\hbar$ such that:
\begin{itemize}
\item
For any $f\in \C^\infty_0(M)$, the section $\hbar \mapsto Q_\hbar(f)$ is continuous.
\item
For each value of $\hbar$, the image of $Q_\hbar$ densely generates $\A_\hbar$.
\item
The map $Q_\hbar$ intertwines complex conjugation (of functions) with the involution (in $\A_\hbar$).
\item
$Q_0$ extends to an isomorphism $Q_0 : \C_0(M)\isom \A_0$.
\item
For any $f,g\in\C^\infty_0(M)$, as $\hbar\to0$ the norm
\[
\left\Vert\tfrac1{i\hbar}[Q_\hbar(f),Q_\hbar(g)] - Q_\hbar(\{f,g\})\right\Vert
\]
converges to $0$.
\end{itemize}
This is not intended as a complete definition, because there are many variations on this in the literature (see \cite{haw10} for a review of the variations). 

Based on the idea of a smoothly deformed product, it was initially assumed that the set of $\hbar$-values would be an interval, and each $Q_\hbar$ would be injective.
Those assumptions are overly restrictive. There is a nice quantization of the symplectic manifold $S^2$ in which the algebras $\A_{\hbar>0}$ are finite-dimensional. In that case, the quantization maps cannot possibly be injective and the dimension cannot vary continuously over an interval of $\hbar$-values. Instead, we should only assume that the set of $\hbar$-values is dense at $0$.

In my view, the main objective of quantization is a single algebra. In these terms, that would be $\A_1$. Strict deformation is flexible enough that $\A_1$ is in no way uniquely determined by the Poisson manifold. However, there may exist a systematic natural construction if we include additional structure.

Imagine quantization as a contravariant functor whose codomain is the category of \cs-algebras and $*$-homomorphisms. The objects of the domain are Poisson manifolds with additional structure and the morphisms are Poisson maps compatible with this additional structure. 

In this picture, the continuous field is not the direct result of quantizing $M$. Instead, the algebra of continuous sections is the quantization of a larger Poisson manifold $M\times \R$ (with the Poisson structure rescaled by $\hbar\in\R$). The evaluation homomorphism from continuous sections to $\A_\hbar$ should be the quantization of the inclusion map $M\times \{\hbar\} \into M\times\R$. In this way, the quantization functor should contain all structure necessary to produce the continuous field for strict deformation quantization.

In this paper, I am only sketching the construction of the \cs-algebra and 
will not be addressing morphisms at all. The point of this is to identify the ``additional structure'' needed for quantization. I am proposing that an object in the domain of quantization consists of a Poisson manifold $M$, a symplectic groupoid $\Sigma$  integrating $M$, a prequantization of $\Sigma$, a polarization of $\Sigma$, and a choice of half-forms. All of these things will be defined below.

The reader may be wondering where the quantization maps are supposed to come from in this picture. 
In strict deformation quantization, the existence of quantization maps is important, but the specific maps are not. If the $Q_\hbar$'s are changed by something of order $\hbar^2$, then this does not effect the algebraic structure.
What the quantization maps really do is to give a little bit of differentiable structure to the continuous field of \cs-algebras. For $f\in \C^\infty_0(M)$, the section $\hbar\mapsto Q_\hbar(f)$ should be seen as differentiable at $\hbar=0$.
Specifying a differentiable (or smooth) subalgebra of continuous sections sections should be almost as good as defining quantization maps.

In practice there are many examples of quantization maps, but their construction involves structures that are specific to classes of examples (such as \Kahler\ manifolds or linear Poisson structures). For these reasons, I do not expect quantization maps to be as natural or general as the algebras. 
In the approach that I am proposing, the \cs-algebra is constructed using a symplectic groupoid. With such a geometrical construction, it is very plausible that there will be a nice ``smooth'' subalgebra of the \cs-algebra; this would compensate for the lack of a general quantization map.

In the end, strict deformation quantization may not be essential to quantization. It provides some of the motivation, but if there is a functor that produces strict deformation quantizations in many but not all cases, then that is still interesting.

I will discuss in Section \ref{Completion} the likely ambiguity in the choice of \cs-completion. This means that there may be more than one quantization functor, depending upon the recipe for choosing a norm. Perhaps only one of these choices will lead to strict deformation quantizations, but other choices may have other interesting properties.

\subsection{Synopsis}
I begin in Section~\ref{Geometric} by reviewing the standard geometric quantization of symplectic manifolds. This constructs a Hilbert space, but the relevant \cs-algebra is that of compact operators. Of course, a symplectic manifold is a special type of Poisson manifold. The algebra can be constructed directly using the pair groupoid, and this case motivates my general approach.

In Section~\ref{Ingredients}, I discuss the ingredients of my quantization recipe. I review the results of Weinstein, Xu, Crainic, and Zhu on the prequantization of symplectic groupoids and fill in some additional details. I also present some related results which will be useful later. Then I present my definition of a groupoid polarization, Definition~\ref{polarization}.

In Section \ref{Algebra}, I discuss the twisted, polarized \cs-algebra itself. I give a preliminary definition which only applies to well-behaved cases; nevertheless, this definition is sufficient for most of the examples. I also extend the definition a little further using the idea of Bohr-Sommerfeld quantization conditions, which are quite standard in geometric quantization. 
A more complete definition is deferred to a future paper.

In Section~\ref{Examples}, I present a series of examples in which my procedure reproduces known examples of quantization of Poisson manifolds.

In Section~\ref{Polarization2}, I present some tools for working with polarizations of groupoids. In particular, for real polarizations of a groupoid, I present the corresponding structure of a real polarization of a Lie algebroid. This leads to a definition for a real polarization of a Poisson manifold. It also reveals some limitations on real polarizations.

Lest the reader suspect that this procedure \emph{only} reproduces known examples, I present further examples in Section~\ref{Examples2}. Some of these are just examples of polarizations, in other cases I do sketch the construction of the algebra.

\section{Groupoid Preliminaries}
I will assume that the reader is familiar with the basic definitions of Lie groupoids. For this background, see \cite{mac3,cds-w,m-m,lan2}. I use the term ``Lie groupoid'', but the terms ``smooth groupoid'' and ``differential groupoid'' are also used synonymously in the literature. All groupoids here are Lie groupoids.

A Lie groupoid \emph{homomorphism} is a smooth map between Lie groupoids that is a functor if the groupoids are viewed as categories.

Given a groupoid $\G$, I denote by $\G_k\subset\G^k$ the $k$-nerve; that is, the set of composable chains of $k$ elements. The $0$-nerve $\G_0$ is the base manifold of $\G$, which I usually call $M$. I use the following notations for the structure maps of a generic groupoid. The unit is $\unit: M\into \G$, but I usually treat $M$ as a submanifold, unless I need to refer to the inclusion explicitly. The source and target maps are $\s,\t:\G\onto M$. The multiplication map is $\m:\G_2\onto\G$, but for a composable pair $(\gamma,\eta)\in\G_2$ I denote multiplication by apposition or a dot as $\gamma\eta=\gamma\cdot\eta=\m(\gamma,\eta)$. The inverse map is $\inv:\G\to\G$, but for $\gamma\in\G$ I denote its inverse as $\gamma^{-1}=\inv(\gamma)$. I also denote the Cartesian projections to the first and second factors of $\G^2$ as $\pr_1,\pr_2:\G_2\subset \G^2\onto\G$.

The nerve of a groupoid has a natural structure as a simplicial manifold. The maps $\s$, $\t$, $\m$, $\pr_1$, and $\pr_2$ are face maps.  This structure gives a simplicial coboundary operator $\partial^*$ on differential forms on the nerve, $\Omega^\bullet(\G_\bullet)$. In particular, for $\theta\in\Omega^\bullet(M)$,
\[
\partial^*\theta := \t^*\theta - \s^*\theta \in \Omega^\bullet(\G)
\]
and for $\omega\in\Omega^\bullet(\G)$,
\begin{subequations}
\beq
\partial^*\omega := \pr_1^*\omega - \m^*\omega + \pr_2^*\omega \in\Omega^\bullet(\G_2) .
\eeq

\begin{definition}
A differential form $\omega\in\Omega^\bullet(\G)$ is \emph{multiplicative} if
\beq
 \partial^*\omega = 0 .
\eeq
A \emph{symplectic groupoid} is a groupoid $\Sigma$ with a multiplicative symplectic form $\omega\in\Omega^2(\Sigma)$. (I usually denote a symplectic groupoid as $\Sigma$.)
\end{definition}
\label{multiplicative.form}
\end{subequations}

It will also be useful to extend the simplicial coboundary to line bundles, in a multiplicative sense. If $\Lambda$ is a line bundle over $M$ and $\G$ is a groupoid over $M$, then
\[
\partial^* \Lambda := \t^*\Lambda \otimes \s^*\Lambda^*
\]
is a line bundle over $\G$. If $L$ is a line bundle over $\G$, then 
\beq
\label{L_coboundary}
\partial^* L := \pr_1^* L \otimes \m^* L^* \otimes \pr_2^* L
\eeq
is a line bundle over $\G_2$. This continues so that the coboundary of a line bundle over $\G_k$ is a line bundle over $\G_{k+1}$. The coboundary of a coboundary, $\partial^*\partial^* L$, is a canonically trivial line bundle. If $L$ is equipped with a connection, then $\partial^*L$ inherits a connection, and 
\[
\curv \partial^*L = \partial^*(\curv L) .
\]
A section of a line bundle $\sigma \in \Gamma(\G_k,L)$ has a multiplicative coboundary $\partial^*\sigma \in \Gamma(\G_{k+1},\partial^* L)$. 

It is standard terminology to use $\s$ as a prefix to describe a property of the fibers of $\s:\G\to M$. Hence $\G$ is \emph{$\s$-connected} if the fibers are connected, \emph{$\s$-simply connected} if they are simply connected, and \emph{$\s$-locally trivial} if they form a bundle. These properties are exactly the same as ``$\t$-connected'' \emph{et cetera}, but it is just conventional to refer to $\s$.

I also use a superscript to denote certain subbundles of the tangent bundle, namely $T^\s\G := \ker T\s \subset T\G$ and likewise $T^\t\G$ and $T^\m\G_2$.

I denote an arbitrary Lie algebroid over $M$ as $A$. Every Lie algebroid anchor map is denoted as $\#:A\to TM$. This includes the anchor map given by a Poisson structure $\#:T^*M\to TM$ and even the map given by the symplectic structure on a symplectic groupoid, $\# : T^*\Sigma \to T\Sigma$. When there is a real danger of confusion, I will use a subscript to indicate which anchor map it is.


I denote by $\Alg$ the Lie functor from the category of Lie groupoids (and smooth homomorphisms) to the category of Lie algebroids, see \cite{mac3}. This restricts to the classical Lie functor from Lie groups to Lie algebras.

I denote by $\Grp$ the inverse of $\Alg$ (see \cite{c-f1}). The domain is the full subcategory of integrable Lie algebroids (i.e., the image of $\Alg$). If $A$ is an integrable Lie algebroid, then $\Grp(A)$ is the unique (up to isomorphism) $\s$-connected and $\s$-simply connected groupoid such that $A\cong \Alg\Grp(A)$. I call $\Grp(A)$ \emph{the} integration of $A$, but I will say that any $\s$-connected groupoid $\G$ with $\Alg(\G)\cong A$ \emph{integrates} $A$.

If $M$ is a Poisson manifold, then the cotangent bundle $T^*M$ with the Koszul bracket is a Lie algebroid. I will say that a symplectic groupoid over $M$ integrates $M$ if $\t$ is a Poisson map; $\Sigma$ is a symplectic integration of $M$ if and only if it is a groupoid integration of $T^*M$.
Any Lie groupoid integrating this has a natural symplectic groupoid structure. I will denote \emph{the} symplectic integration by $\Sig(M):=\Grp(T^*M)$ (see \cite{c-f2}).

Any manifold can be regarded as a groupoid, with every point a unit. Given a manifold, we can also construct the pair groupoid $\Pair(M):=M\times M$; note that this integrates $TM$. If $M$ happens to be a symplectic manifold, then the pair groupoid is a symplectic groupoid if we give it the symplectic structure of $M\times M^-$, where $M^-$ is $M$ with the symplectic form reversed. In this case $\Pair(M)$ is a symplectic integration of $M$.

For an arbitrary manifold, there is also the fundamental groupoid $\PI(M)$. This is the $\s$-simply connected cover of $\Pair(M)$, hence $\PI(M)=\Grp(TM)$. If $M$ is symplectic, then $\PI(M)$ inherits a symplectic groupoid structure from $\Pair(M)$, and $\Sig(M)=\PI(M)$.

Given a Lie groupoid $\G$ over $M$, we can construct another Lie groupoid simply by applying the tangent functor. The unit manifold of $T\G$ is $TM$, the source map is $T\s:T\G\to TM$, the multiplication is $T\m: T\G_2\to T\G$, \emph{et cetera}.  The complexified tangent bundle $T_\co\G$ is a groupoid in the same way. This tangent groupoid is very important for my definition of polarization, but it should not be confused with Connes' tangent groupoid \cite{con1}, which is also used in an example of quantization.

The cotangent bundle $T^*\G$ is a groupoid as well, but in a very different way. This is a symplectic groupoid over $\Alg^*(\G)$, the dual vector bundle of the Lie algebroid. This is an important example in Section~\ref{Linear}.

\section{Geometric Quantization}
\label{Geometric}
For background on geometric quantization, see \cite{a-e,b-w,sni,woo}.
Let $(M,\omega)$ be a symplectic manifold. Geometric quantization constructs a Hilbert space based on $(M,\omega)$. The \cs-algebra is that of compact operators on this Hilbert space. However, the construction requires a little more structure than just the symplectic form.

The first step of the geometric quantization procedure is known as \emph{prequantization}. In the formulation of Kostant, this means the construction of a Hermitian line bundle $L$ over $M$ with connection $\nabla$ and curvature $\omega$. Such an $L$ exists if and only if the cohomology class of $\omega/2\pi$ is integral. There is an equivalent formulation due to Souriau using a circle bundle instead. This is just the principal $\T$-bundle associated to $L$.

\subsection{Polarization}
The space of smooth sections of $L$ is too large for geometric quantization. We can restrict to a smaller class of sections by using a polarization.
\begin{definition}
A \emph{polarization} of a symplectic manifold $(M,\omega)$ is an involutive (i.e., integrable) Lagrangian distribution $\F\subset\TcM$. A section $\psi\in\Gamma(M,L)$ is \emph{polarized} if 
\beq
\label{polarized}
0=\nabla_X\psi \quad\forall X\in\F .
\eeq 
I define polarized sections for other line bundles in the same way.
\end{definition}

For two vector fields $X,Y\in\Gamma(M,\F)$ and a polarized section $\psi\in\Gamma(M,L)$, the definition of curvature shows that
\[
0 = [\nabla_X,\nabla_Y]_-\psi = \nabla_{[X,Y]}\psi + i\omega(X,Y)\psi .
\]
The involutivity of $\F$ implies that the first term vanishes. Since $\F$ is Lagrangian, it is in particular isotropic and this means that the second term vanishes. ``Lagrangian'' is equivalent to ``isotropic of maximal dimension''. We can interpret the definition of polarization to mean that \eqref{polarized} is as restrictive as possible while still being consistent.

\begin{definition}
 If $\F\cap\bar\F=0$ then $\F$ is \emph{totally complex} (or \emph\Kahler). A \emph{real polarization} is an involutive Lagrangian real distribution $\F\subset TM$.
\end{definition}
If a (complex) polarization satisfies $\bar\F = \F$, then it is the complexification of a real polarization. In this way, I will regard real as a special case of complex.

In general, the real part, $\F\cap TM$ is not necessarily of constant rank, but if it is then it is an involutive distribution, i.e., a foliation. In that case, $\F+\bar\F \subseteq T_\co M$ is a subbundle, but is not necessarily involutive. We can choose to assume that these are very well behaved.

\begin{definition}
\label{admissible1}
The (generally singular) distributions $\D,\E\subseteq TM$ are defined by 
$\D_\co :=\F\cap\bar\F$ and $\E_\co:=\F+\bar\F$. A polarization $\F$ is \emph{strongly admissible} if there exist manifolds and surjective submersions $M \stackrel{p}\onto M/\D \stackrel{q}\onto M/\E$ such that $\D$ and $\E$ are the  the kernel foliations, $\D= \ker Tp$ and $\E = \ker T(q\circ p)$.
\end{definition}

\subsection{Hilbert Space}
The next step is to construct a Hilbert space $L^2_
\F(M,L)$ based on the polarized sections of $L$. Assume that $\F$ is a strongly admissible. There are two important subtleties. 

The naive approach would be to define $L^2(M,L)$ using the symplectic volume form, and then to consider the Hilbert subspace densely spanned by polarized sections. However, a polarized section is covariantly constant along the leaves of $\D$, and if these leaves are noncompact, then polarized sections are not square-integrable. This is solved by recognizing that a polarized section is equivalent to a section over $M/\D$. However, there is (in general) no natural volume form on $M/\D$, so we must absorb the volume form into the choice of line bundle.

Before considering this, we should recall a couple of standard concepts which will be important throughout this paper.
\begin{definition}
\label{square_root}
If $\Lambda$ is a (real or complex) line bundle, then \emph{a square root} of $\Lambda$ is another line bundle $\sqrt\Lambda$ equipped with a specific isomorphism $\sqrt\Lambda\otimes\sqrt\Lambda\cong\Lambda$. 
\end{definition}
This may not exist, and if it does \emph{it may not be unique}. For $\Lambda$ a complex line bundle, $\sqrt{\Lambda}$ exists if and only if $c_1(\Lambda) = 0 \in H^2(M;\Z_2)$. For a real line bundle, $c_1(\Lambda)=0$, so $\sqrt{\Lambda_\co}$ always exists, but a real square root $\sqrt{\Lambda}$ only exists if $\Lambda$ is orientable (i.e., trivial); in that case there exists a unique preferred choice of $\sqrt\Lambda$ which is real and orientable. Think of this as the positive square root. For both real and complex line bundles, the nonuniqueness of $\sqrt{\Lambda}$ is the freedom to take a tensor product with a \emph{real} line bundle, and this is classified by $H^1(M;\Z_2)$. By the Leibniz rule, any connection on $\Lambda$ determines an equivalent connection on $\sqrt{\Lambda}$.

\begin{definition}
If $\F\subset T_\co M$ is a distribution, then its \emph{annihilator} (or \emph{conormal}) bundle is 
\[
\F^\perp := \{\xi \in T^*M \mid \forall X\in \F : \langle X,\xi\rangle =0\} .
\]
\label{Bott}%
If $\F$ is an involutive distribution, then the \emph{Bott connection} is the flat $\F$-connection on $\F^\perp$ that equals the Lie derivative. That is,
\[
\nabla_X \xi := \Lie_X \xi = X\inner d\xi
\]
for any $X\in \Gamma(M,\F)$ and $\xi \in \Gamma(M,\F^\perp)$. The Bott connection extends to any bundle constructed from $\F^\perp$, including $\sqrt{\Wedgem\F^\perp}$.
\end{definition}

For a moment, consider the case that $\F$ is a real and strongly admissible polarization. Pulling back $\Wedgem T^*(M/\F)$ to $M$ gives the bundle $\Wedgem \F^\perp$. Now (if possible) construct the positive square root bundle $\sqrt{\Wedgem \F^\perp}$. A section of $\sqrt{\Wedgem \F^\perp}$ is the pull-back of a section of $\sqrt{\Wedgem T^*(M/\F)}$ if and only if it is $\F$-constant by the Bott connection, i.e., it is polarized.

Combining the connection on $L$ with the Bott connection gives a flat $\F$-con\-nec\-tion on $L\otimes\sqrt{\Wedgem \F^\perp}$. Define a section of $L\otimes\sqrt{\Wedgem \F^\perp}$ to be \emph{polarized} if it is $\F$-constant by this connection. The local inner product of two polarized sections is an $\F$-constant section of $\Wedgem \F^\perp_\co$; this can be integrated over over $M/\F$ to define the global inner product for $L^2_\F(M,L)$.

If $\F$ is not real, then the local inner product is valued in $\sqrt{\Wedgem \bar \F^\perp}\otimes \sqrt{\Wedgem \F^\perp}$ rather than $\Wedgem\D^\perp_\co$. However, there is a natural isomorphism 
\[
\Wedgem \bar \F^\perp \otimes \Wedgem \F^\perp \cong \Wedgem\E^\perp_\co \otimes \Wedgem\D^\perp_\co \mbox,
\]
and the exterior product with $\frac{\omega^k}{k!}$ (where $2k:=\rk \E - \rk \D$) defines a canonical isomorphism $\Wedgem\E^\perp \cong \Wedgem\D^\perp$. This corrects the inner product and gives a (fairly) general definition of $L^2_\F(M,L)$. 

The bundle $\sqrt{\Wedgem\F^\perp}$ is known as a ``half-form bundle''. Although I have not indicated it explicitly, $L^2_\F(M,L)$ does depend on this. The subscript $\F$ should be understood to represent the choice of polarization \emph{and} half-forms.

The second subtlety is that a polarized section of $L\otimes\sqrt{\Wedgem \F^\perp}$ is covariantly constant along $\D$-leaves. This is locally consistent, because $\D$ is isotropic and hence the connection restricted to a $\D$-leaf is (locally) flat. However, if a $\D$ leaf is not simply connected, then the connection determines a holonomy homomorphism $\pi_1(\text{leaf})\to\mathbb{T}$. If this is nontrivial, then any polarized section must vanish on this leaf.
 
When there is nontrivial holonomy, it varies continuously over $M/\D$ and is typically nontrivial over a dense subset. This means that global polarized sections must vanish over a dense subset. Continuity implies triviality. 

The usual solution to this paradox is to instead use \emph{distributional} polarized sections supported on the union of $\D$-leaves with trivial holonomy, a subset of $M$ known as the \emph{Bohr-Sommerfeld} variety, $M_{\BS}$. In practice, this means working with sections over $M_{\BS}/\D$.

This is a pragmatic solution, but it lacks a theoretical justification and it doesn't always work. There are other reasons that the Hilbert space of polarized sections may be ``too small''. For example, if the leaves of a real polarization are dense or a \Kahler\ metric is negative definite, then there may be no global polarized sections.

There is a plausible solution to all of this using so called ``cohomological wave functions'' (see \cite{woo}). Local polarized sections of $L\otimes\sqrt{\Wedgem\F^\perp}$ form a sheaf. The space of global polarized sections is the $0$'th sheaf cohomology space. It has been proposed that the correct Hilbert space is a completion of the \emph{total} sheaf cohomology. The only problem is to construct a natural inner product and this has not been solved.

However, we can conclude that the simple construction of $L^2_\F(M,L)$ from the smooth, polarized, global sections is probably correct if and only if the higher degree cohomology vanishes.

\begin{example}
The canonical example of a symplectic manifold is a cotangent bundle $T^*N$. This is the phase space for a physical system with configuration space $N$. The canonical choice of polarization $\F$ is the foliation whose leaves are the cotangent fibers. In this case there are no Bohr-Sommerfeld conditions and the prequantization $L$ is trivializable, so the Hilbert space is
\[
L^2_\F(T^*N,L) \cong L^2_\F(T^*N) \cong L^2(T^*N/\F) \cong L^2(N) .
\]
The polarized phase space behaves like the configuration space.
\end{example}

This is the simplest instance of an important guiding principle. In the simplest cases, a polarized space behaves like its quotient by the polarization.

\subsection{Algebra}
\label{GQ.Algebra}
As I stated in the introduction, the quantization problem that concerns me here is the construction of a noncommutative \cs-algebra that is approximated by the Poisson algebra of functions on $M$. 
For standard geometric quantization of a symplectic manifold, that algebra is the \cs-algebra of compact operators, $\K[L^2_\F(M,L)]$. 

This is very specific to symplectic manifolds. A generic Poisson manifold will not quantize to the algebra of compact operators on a Hilbert space. In order to generalize geometric quantization from the symplectic case, we need to break free of the Hilbert space and construct the algebra directly. This can be done by working over the (symplectic) pair groupoid $\Pair(M)=M\times M^-$. 

Let's temporarily restrict to the well-behaved case, so that the Hilbert space $L^2_\F(M,L)$ is densely spanned by smooth, polarized sections of $L\otimes\sqrt{\Wedgem\F^\perp}$. The algebra $\K[L^2_\F(M,L)]$ is densely spanned by ``ket-bras'' --- products over $\Pair(M)$ of a section and a complex conjugate section. There is a dense subalgebra whose elements are given by polarized sections of 
\[
\left(L\otimes \sqrt{\Wedgem\F^\perp}\right) \boxtimes \left(\bar L \otimes \sqrt{\Wedgem\bar\F^\perp}\right)
\]
over $\Pair(M)$. The $\boxtimes$ denotes the outer tensor product of vector bundles over a Cartesian product space.

This algebra is a \emph{twisted} and \emph{polarized} convolution algebra over $\Pair(M)$. It relates to the convolution algebra in the same way that $L^2_\F(M,L)$ relates to $L^2(M)$.  

The ``twist'' is the line bundle $L\boxtimes\bar L$. Its curvature is the symplectic form of $\Pair(M)$. So, $L\boxtimes\bar L$ serves as a prequantization of the symplectic groupoid $\Pair(M)$.

This bundle has more structure than just its connection. The inner product on $L$ allows us to multiply vectors in different fibers of $L\boxtimes\bar L$; denote this product as $\sigma$. This is an elementary example of a groupoid cocycle. 
In general, a twisted groupoid \cs-algebra is constructed using a line bundle and a cocycle; see \cite{ren} and Section~\ref{Twisted}.

The second modification is polarization. The modified convolution algebra of the pair groupoid is composed of sections which are polarized by 
\[
\F\times\bar\F \subset T_\co(\Pair M) .
\] 
I will denote this algebra as
\[
\K[L^2_\F(M,L)] = \cs_{\F\times\bar\F}(\Pair M,\sigma) .
\]

This suggests my general recipe for quantization of a Poisson manifold $M$. The algebra should be a twisted, polarized groupoid algebra
\[
\cs_\P(\Sigma,\sigma)
\]
where $\Sigma$ is a symplectic groupoid integrating $M$, $\P$ is a polarization of $\Sigma$, and $\sigma$ denotes a prequantization of $\Sigma$.

All of the issues that must be considered in constructing the inner product over $M$ translate into issues in the construction of the convolution product over the groupoid. 
Because of the polarization, we do not have the luxury of starting with compactly supported sections. Likewise, if $\D\neq 0$, then the domain of integration in the definition of convolution must be modified. The holonomy problem remains as well and may force us to start with distributional sections.

\section{The Ingredients}
\label{Ingredients}
Before considering the recipe for the algebra, we need to prepare the ingredients.
\subsection{Integration}
The first ingredient is a symplectic groupoid.

If $\Sigma$ is a symplectic groupoid over a manifold $M$, then there exists a unique Poisson structure on $M$ such that $\t : \Sigma \to M$ is a Poisson map. The Lie algebroid $\Alg(\Sigma)$ is naturally identified with the cotangent bundle $T^*M$; the anchor map $\#:T^*M\to TM$ is given by contraction with the Poisson bivector $\pi$, and the Lie algebroid bracket is the Koszul bracket of $1$-forms. 

We can look at this from the other direction, by starting with a given Poisson manifold $M$ and trying to construct a symplectic groupoid over it. The Lie algebroid structure on $T^*M$ is fixed by the Poisson structure, so the groupoid structure is given locally by the Poisson structure. If the $\s$-fibers of a symplectic groupoid are not connected, then there is some global structure which is not given by the Poisson structure. That is not relevant to the quantization of $M$, so I will assume that the groupoid is $\s$-connected.
\begin{definition}
If $M$ is a Poisson manifold, then a symplectic groupoid $\Sigma$ over $M$ \emph{integrates} $M$ if $\t:\Sigma\to M$ is a Poisson map and $\Sigma$ is $\s$-connected. A Poisson manifold $M$ is \emph{integrable} if such a symplectic groupoid exists.
\end{definition}
A complete understanding of integrability was achieved by Crainic and Fernandes \cite{c-f2}. They presented a computable necessary and sufficient condition for integrability. A Poisson manifold is integrable if and only if $T^*M$ is integrable to a Lie groupoid. If $\Sigma$ integrates a Poisson manifold $M$, then $\t:\Sigma\to M$ is a complete symplectic realization, and conversely, if a complete symplectic realization of $M$ exists, then $M$ is integrable.

Unfortunately, this notion of integrability does not include the usual cozy assumption of Hausdorffness. In the theory of Lie groupoids, it is often necessary to consider non-Hausdorff groupoids. However, the base manifolds, algebroids, and the $\s$ and $\t$ fibers are usually assumed to be Hausdorff. Non-Hausdorff groupoids raise many delicate issues, which I do not attempt to address in this paper.

If $M$ is an integrable Poisson manifold, then there exists a \emph{unique} $\s$-simply connected symplectic groupoid integrating $M$.
\begin{definition}
If $M$ is an integrable Poisson manifold, then \emph{the symplectic integration} $\Sig(M)$ is the unique $\s$-simply connected symplectic groupoid integrating $M$.
\end{definition}
As a groupoid, this is $\Sig(M)\cong\Grp(T^*M)$ the unique $\s$-simply connected groupoid integrating the Lie algebroid $T^*M$. Any symplectic groupoid integrating $M$ is a quotient of $\Sig(M)$.

It would be tempting to just use $\Sig(M)$ in the quantization construction. However, in the motivating example of a symplectic manifold $M$, the relevant symplectic groupoid was $\Pair(M)$. This is a symplectic groupoid integrating $M$, but it is not isomorphic to $\Sig(M)$ unless $M$ is simply connected.

So, to quantize a Poisson manifold $M$, the first ingredient in this recipe is \emph{some} symplectic groupoid $\Sigma$ integrating $M$.

This entails existence and uniqueness issues. Such a groupoid only exists if $M$ is integrable. This should be seen as a potential obstruction to quantization. In general, the symplectic groupoid is not unique, although all possible choices are quotients of $\Sig(M)$. This should be viewed as an ambiguity in the quantization process.

\subsection{Prequantization}
\label{Prequantization}
Prequantization of symplectic groupoids has been studied by Weinstein and Xu \cite{w-x}, by Crainic \cite{cra2}, and by Crainic and Zhu \cite{c-z}. I summarize their relevant results and fill in some more details. 

Prequantizing a symplectic groupoid $\Sigma$ involves a little more structure than simply prequantizing $\Sigma$ as a symplectic manifold. There are again two different but equivalent perspectives on prequantization: the Souriau picture (in terms of circle bundles) and the Kostant picture (in terms of line bundles). The contrast between these perspectives becomes greater, but both perspectives are useful.
I favor the Kostant picture, because it is better suited to the construction of the algebra. 

The purpose of prequantizing a symplectic groupoid is eventually to construct twisted polarized convolution algebras, but it is useful to first concentrate on the structure needed just to define a twisted convolution algebra. This is meaningful for an arbitrary groupoid.

Let $L$ be a line bundle over a groupoid $\G$. In order to define a convolution algebra with coefficients in $L$, we need an associative way of multiplying fibers of $L$ over different points; this means, for a composable pair of groupoid elements $(\gamma,\eta)\in\G_2$, a bilinear map $\sigma(\gamma,\eta): L_\gamma\otimes L_\eta \to L_{\gamma\eta}$. To construct a \cs-algebra, we also need a norm and an (antilinear) adjoint ${}^*:L_\gamma \to L_{\gamma^{-1}}$. A vector bundle with these structures is known as a \emph{Fell bundle} \cite{fel,kum3,yam}. 

We can equivalently think of $\sigma$ as a section of 
\beq
\partial^*L^* = \pr_1^*L^*\otimes \m^*L\otimes\pr_2^*L^*
\eeq
over $\G_2$, see \eqref{L_coboundary}. Associativity means that the (multiplicative) coboundary of $\sigma$ equals $1$.

Because $L$ is a line bundle, the norm and adjoint are equivalent to a Hermitian inner product on $L$. The cocycle $\sigma$ must have norm $1$ everywhere. So the structure we need can be summarized as: a Hermitian line bundle with a norm $1$ cocycle. This is the Kostant picture of a twist.

There is an equivalent Souriau picture. A Hermitian line bundle is equivalent to a circle bundle (the set of elements in $L^*$ with norm $1$). A cocycle $\sigma$ gives this circle bundle the structure of a groupoid, which I denote as $\G^\sigma$. In fact, it is a $\T$-extension:
\beq
\label{extension}
\T\times M \into\G^\sigma \onto \G .
\eeq
This is a short exact sequence of groupoids in the following sense: As a subgroupoid, $\T\times M$ acts on $\G^\sigma$ (from the left, say); the second map is a fibration (see Def.~\ref{fibration} below) and its fibers are the orbits of the $\T\times M$ action. Any $\T$-extension is given by a Hermitian line bundle and multiplicative cocycle.
The set of isomorphism classes of $\T$-extensions of $\G$ forms a group $\Tw(\G)$; see \cite{kum1,kum2,k-m-n-w}.

Given an extension, $\G^\sigma$, applying the functor $\Alg$ to \eqref{extension} gives a Lie algebroid extension,
\beq
\label{algebroid_extension}
0 \to \R\times M \to \Alg(\G^\sigma) \to \Alg(\G) \to 0 .
\eeq
As a short exact sequence of vector bundles, this can be split, and a splitting identifies $\Alg(\G^\sigma) \cong \Alg(\G)\oplus (\R\times M)$ as vector bundles. 
With such an identification, the bracket is of the form,
\beq
\label{ext.bracket}
[(\xi,f),(\zeta,g)]_{\Alg{\G^\sigma}} = ([\xi,\zeta]_{\Alg\G},\Lie_{\#\xi}g-\Lie_{\#\zeta}f + c(\xi,\zeta)) .
\eeq
The term $c$ is a Lie algebroid cocycle.
This defines a characteristic class map $\Psi:\Tw(\G) \to H^2_{\mathrm{Lie}}(\Alg\G), [\sigma]\mapsto [c]$ to Lie algebroid cohomology; see \cite{cra2}.

\begin{thm}
\label{Psi_injective}
For $\G$ any $\s$-simply connected Lie groupoid,  $\Psi:\Tw(\G) \to H^2_{\mathrm{Lie}}(\Alg\G)$ is injective.
\end{thm}
\begin{proof}
Let $\G^\sigma$ be some $\T$-extension of $\G$. The class $\Psi(\sigma)$ determines the Lie algebroid extension \eqref{algebroid_extension}
modulo isomorphism. Obviously, all three terms are integrable Lie algebroids, so we can apply the integration functor $\Grp$ to this Lie algebroid extension and get a sequence of groupoids,
\beq
\label{integrated_extension}
\R\times M \to \Grp\Alg(\G^\sigma) \onto \G .
\eeq
The first map is \emph{not} necessarily injective. The groupoid $\Grp\Alg(\G^\sigma)$ is the $\s$-simply connected cover of $\G^\sigma$. In order to recover a $\T$-extension from this, quotient $\R$ by $2\pi\,\Z$ and $\Grp\Alg(\G^\sigma)$ by the image of $2\pi\,\Z$. The result is isomorphic to the extension~\eqref{extension} for $\G^\sigma$.
\end{proof}

Now we can turn to prequantization of symplectic groupoids. A prequantization of $\Sigma$ should be both a prequantization of $\Sigma$ as a symplectic manifold and a twist of $\Sigma$ as a groupoid. In order for convolution to be compatible with polarization, the twist must be compatible with the connection on $L$; this means that the cocycle $\sigma$ should be a covariantly constant section of $\partial^*L$ over $\Sigma_2$.

\begin{definition}
A \emph{prequantization of a symplectic groupoid} $\Sigma$ consists of a Hermitian line bundle $L\to\Sigma$ with connection and a section $\sigma\in\Gamma(\Sigma_2,\partial^*L^*)$ such that:
\begin{enumerate}
\item
The curvature of $L$ equals the symplectic form;
\item
$\sigma$ is a (multiplicative) cocycle and has norm $1$ at every point; 
\item
$\sigma$ is covariantly constant.
\end{enumerate}
\end{definition}
As a shorthand, I will usually refer to a prequantization $(\sigma,L,\nabla)$ simply 
as $\sigma$. This avoids clutter in the notation for a twisted polarized convolution \cs-algebra.

In the Souriau picture, a prequantization of $\G$ is a $\T$-extension of $\G$ which (as a $\T$-bundle) is equipped with a connection that has curvature $\omega$ and is compatible with the groupoid structure.

The curvature of $\partial^*L^*$ is $-\partial^*\omega =0$, because the symplectic form is multiplicative. A covariantly constant cochain $\sigma\in\Gamma(\Sigma_2,\partial^*L^*)$ therefore exists if and only if the holonomy of $\partial^*L^*$ is trivial. If $\Sigma_2$ is connected, then $\sigma$ is unique up to an irrelevant multiplicative constant.

Consider $\unit^*L$, the restriction of $L$ to the identity submanifold $M\subset \Sigma$. The cocycle $\sigma$ makes $\unit^*L$ into a bundle of algebras, each isomorphic to $\co$. There is thus a unit section of $\unit^*L$ over $M$. Compatibility of $\sigma$ with the connection implies that this unit section is covariantly constant, so the connection $\unit^*L$ must be flat and have trivial holonomy. In fact, one of the elementary properties of a symplectic groupoid is that the unit submanifold is Lagrangian, so the curvature of $\unit^*L$ is automatically
\[
\curv \unit^*L = \unit^* \omega = 0 .
\]
So, the only necessary condition for this unit section to exist is that $\unit^*L$ have trivial holonomy.

Weinstein and Xu showed that this condition is not just necessary for prequantization --- it is sufficient.
\begin{thm} 
\label{WX_prequantization}
\cite{w-x} If $\Sigma$ is a symplectic groupoid and $L\to \Sigma$ is a Hermitian line bundle with a connection such that the curvature equals the symplectic form and $\unit^*L$ has trivial holonomy, then any covariantly constant norm $1$ section of $\unit^*L$ is the unit section for  a unique cocycle making this a prequantization of $\Sigma$.
\end{thm}
They also proved that if $\Sigma$ is $\s$-connected, $\s$-simply connected, $\s$-locally trivial, and (symplectically) prequantizable, then such an $L$ exists and is unique. 

Crainic \cite{cra2} proved a more general prequantizability result by studying the integrability of the Lie algebroid extension \eqref{algebroid_extension}. This was further improved by Crainic and Zhu \cite{c-z}. See also \cite{b-c-z}. Suppose that $\varphi : S^2 \to M$ is a smooth map whose image lies in a single symplectic leaf of a Poisson manifold $M$. That symplectic leaf has a well defined symplectic form $\omega_{\mathrm{Leaf}}$ which can be pulled back to $S^2$. The monodromy \cite{c-f2} of $\varphi$ is the first order variation of the integral $\int_{S^2}\varphi^*\omega_{\mathrm{Leaf}}$. The map $\varphi$ has trivial monodromy if this integral is unchanged to first order if $\varphi$ is perturbed. Equivalently, $\varphi$ has trivial monodromy if it is the base map of a homomorphism of Lie algebroids $TS^2\to T^*M$.
\begin{definition}
\cite{c-z} The \emph{periods} of a Poisson manifold $M$ are the integrals 
\[
\int_{S^2}\varphi^*\omega_{\mathrm{Leaf}}
\]
for such smooth maps $\varphi :S^2\to M$ with trivial monodromy.
\end{definition}

\begin{thm}
\label{c_prequantizability}
\cite{cra2,c-z}
For an integrable Poisson manifold $M$, the symplectic groupoid $\Sig(M)$ is prequantizable if and only if all the periods of $M$ are integer multiples of $2\pi$. If so, the prequantization of $\Sig(M)$ is unique.
\end{thm}

If a symplectic groupoid is not $\s$-simply connected, then prequantization may not be unique. This nonuniqueness is described as follows.
\begin{thm}
\label{nonunique}
Let $\Sigma$ be a prequantizable symplectic groupoid over $M$. Any prequantization of $\Sigma$ determines a bijection from $H^1(\Sigma,M;\T)$ to the set of isomorphism classes of prequantizations of $\Sigma$.
\end{thm}
\begin{proof}
By Theorem \ref{WX_prequantization}, a prequantization of $\Sigma$ is equivalent to a Hermitian line bundle with curvature $\omega$, trivial holonomy over $M$, and a choice of covariantly constant unit section of $\unit^*L$.

Choosing a unit section is equivalent to choosing a value in each connected component of $M$. If the line bundles for two prequantizations are isomorphic, then there is enough freedom in the choice of isomorphism to identify the unit sections. So, the isomorphism classes of prequantizations are in bijective correspondence with the isomorphism classes of Hermitian line bundles with curvature $\omega$ and trivial holonomy over $M$.

If $L$ and $L'$ are two such bundles, then their ratio $L'\otimes L^*$ is a flat Hermitian line bundle with trivial holonomy over $M$. Conversely, if $\Lambda$ is a flat Hermitian line bundle with trivial holonomy over $M$, then the tensor product $\Lambda\otimes L$ is another Hermitian line bundle with curvature $\omega$ and trivial holonomy over $M$. So, picking one prequantization gives a bijective correspondence between isomorphism classes of prequantizations and such flat line bundles.
The set of (isomorphism classes of) flat Hermitian line bundles over $\Sigma$ is equivalent to the cohomology group $H^1(\Sigma;\T)$. 

The target map $\t:\Sigma\to M$ is a right inverse of the unit map $\unit$, so $M$ is a retract of $\Sigma$. This means that $\unit^* : H^\bullet(\Sigma;\T)\onto H^\bullet(M;\T)$ is surjective and the cohomology long exact sequence breaks into short exact sequences, including
\[
0 \to H^1(\Sigma,M;\T) \to H^1(\Sigma;\T) \stackrel{\unit^*}\to H^1(M;\T) \to 0 .
\]

The flat bundles we are interested in are those which become trivial when restricted to $M$. These are classified by the kernel of $\unit^*:H^1(\Sigma;\T) \to H^1(M;\T)$, but the exact sequence shows that this is just the relative cohomology group $H^1(\Sigma,M;\T)$.
\end{proof}
Note that if $\Sigma$ is $\s$-connected and $\s$-locally trivial, then this relative cohomology is the cohomology of the $\s$-fiber.

For a Poisson manifold, the Lie algebroid cohomology of $T^*M$ is canonically isomorphic to the Poisson cohomology of $M$, $H^\bullet_{\mathrm{Lie}}(T^*M)= H^\bullet_\pi(M)$ \cite[Lem.~2.1]{w-x}. The characteristic class of a prequantization of a symplectic groupoid is just the class of the Poisson bivector itself $[\pi]\in H^2_\pi(M)$. In this way, the problem of prequantization is partly a matter of finding an element of $\Tw(\G)$ in the preimage $\Psi^{-1}[\pi]$.

\subsection{Polarization of Groupoids}
Let $\G$ be an arbitrary Lie groupoid, and $\P\subset T_\co\G$ an involutive distribution. My objective is to construct a \cs-algebra from polarized (covariantly $\P$-constant) sections of some line bundle over $\G$. 

There are many technical issues in constructing this algebra, but I will completely ignore these for a moment in order to derive the necessary conditions on $\P$ for such an algebra to be plausible. 
In this spirit, I am temporarily pretending that the product is simply convolution on $\G$; the product is actually more complicated, but that discussion is postponed until Section~\ref{Algebra}.

Suppose that $a$ and $b$ are both polarized sections, that is $\forall X\in\P$: $0=\nabla_Xa=\nabla_X b$. A \cs-algebra should be closed under multiplication, so the first question is, is the convolution $a*b$ polarized? 

The convolution product $a*b$ is defined by integrating $\pr_1^*a\pr_2^*b$ over the fibers of $\m:\G_2\to\G$.
Taking the derivative $\nabla_X(a*b)$ corresponds to taking the derivative of $\pr_1^*a\pr_2^*b$ by some (any) vector $\tilde X\in T_\co\G_2$ such that $T\m(\tilde X)=X$. The integration can absorb differentiation along $T^\m_\co\G_2:=\ker T\m$, and this is reflected in the freedom to choose $\tilde X$ modulo $T^\m_\co\G_2$. Now the question becomes,
\[
0 \stackrel?= \nabla_{\tilde X} (\pr_1^*a\pr_2^*b) = \pr_1^*(\nabla_{T\pr_1(\tilde X)}a)\pr_2^* b + \pr_1^*a\,\pr_2^* (\nabla_{T\pr_2(\tilde X)}b) .
\]
This is satisfied if $T\pr_1(\tilde X),T\pr_2(\tilde X)\in\P$, so $\tilde X$ expresses $X\in\P$ as a product (in the tangent groupoid $T_\co\G$) of two other vectors from $\P$.

On the other hand, the multiplication map for a \cs-algebra is very nearly surjective. (Any element is a finite sum of products.) This suggests that if $0=\nabla_X(a*b)$ for any polarized $a$ and $b$, then $X\in\P$. This means that any product of vectors from $\P$ is also in $\P$; in other words, $\P$ is multiplicatively closed.

Putting these two conditions together gives the following definition for compatibility between $\P$ and the groupoid product.
\begin{definition}
For any distribution $\P\subset T_\co\G$, denote $\P_2:= (\P\times\P)\cap T_\co\G_2$. $\P$  is \emph{multiplicative} if for any $(\gamma,\eta)\in\G_2$:
\[
T\m \left(\P_{2\,(\gamma,\eta)}\right) = \P_{\gamma\eta} \subset T^\co_{\gamma\eta}\G .
\]
\end{definition}

Tang \cite{tan} has also given a definition for a ``multiplicative distribution'' over a groupoid. This is formulated in terms of paths tangent to $\P$, so it only applies to real distributions. For real distributions, his definition is equivalent to half of my definition: the condition that $T\m \left(\P_{2\,(\gamma,\eta)}\right) \supseteq \P_{\gamma\eta}$.

\begin{example}
A vector field $X\in\X^1(\G)$ is multiplicative if it is a groupoid homomorphism when viewed as a map $X:\G \to T_\co\G$. A nonvanishing vector field spans a (rank $1$) distribution. If the vector field is multiplicative, then this distribution is multiplicative.
\end{example}
\begin{example}
For any Lie groupoid $\G$, the source and target tangent bundles, $T^\s\G$ and $T^\t\G$ are multiplicative distributions.
\end{example}
 
Note that this definition can also be expressed (more simply) in terms of the annihilator $\P^\perp\subset T^*_\co\G$. $\P$ is multiplicative if and only if
\[
\m^*\P^\perp = \pr_1^*\P^\perp + \pr_2^*\P^\perp .
\]
This is obviously similar to the condition \eqref{multiplicative.form} for a differential form to be multiplicative.
\begin{example}
For any nonvanishing $1$-form $\theta\in\Omega^1(\G)$, the set of vectors normal to $\theta$ is a distribution. If $\theta$ is a multiplicative $1$-form, then this normal distribution is multiplicative.
\end{example}

\begin{example}
In the Souriau picture, prequantization of a symplectic groupoid $\Sigma$ gives a $\T$-extension $\Sigma^\sigma$ with a compatible connection. The connection can be expressed as a horizontal distribution on $\Sigma^\sigma$. The compatibility between the connection and groupoid structure means precisely that this distribution is multiplicative. This is an example of a multiplicative distribution which is \emph{not} involutive.
\end{example}

Multiplicativity is the compatibility of a distribution with the groupoid multiplication. In order to construct a $*$-algebra, we also need compatibility with the groupoid inverse. In a convolution algebra, the adjoint is $a^*:=\inv^*\bar a$, where $\bar a$ is the complex conjugate, and $\inv:\G\to\G$ is the groupoid inverse map. If $a$ is polarized, then $a^*$ should also be polarized: $\forall X\in \P$
\[
0\stackrel?= \nabla_X a^* = \nabla_X\inv^*\bar a = \inv^*\overline{\nabla_{T\inv\bar X}a} .
\]
This is true if $T\inv (X)\in \bar \P$. So, we should require: 
\begin{definition}
A distribution $\P\subset T_\co\G$ is \emph{Hermitian} if $T\inv(\P)=\bar \P$.
\end{definition}

This leads to my main definition.
\begin{definition}
\label{polarization}
A \emph{polarization of a Lie groupoid} $\G$ is an involutive, multiplicative, Hermitian distribution $\P\subset T_\co\G$.  A \emph{polarization of a symplectic groupoid} $\Sigma$ is a polarization in both the symplectic and groupoid senses --- that is, an involutive, multiplicative, Hermitian, Lagrangian distribution. 
\end{definition}

The issues of existence and uniqueness for groupoid polarizations are largely unexplored. I begin to investigate these questions in Section~\ref{Polarization2}, but this is mostly left to future work.

\subsection{Real and Strongly Admissible Polarizations}
\label{Real1}
\begin{definition}
A \emph{real polarization} (of a symplectic manifold, groupoid, or symplectic groupoid) is a real distribution $\P$ whose complexification $\P_\co$ is a polarization.
\end{definition}
The complexification $\P_\co$ is preserved by complex conjugation, and any complex distribution that is preserved by complex conjugation is the complexification of a real distribution. In this way, I consider real polarizations to be a special case of complex polarizations. 

\begin{definition}
\label{LA-groupoid}
A \emph{Lie algebroid-groupoid} (defined in \cite{mac1}, but called an $\mathcal{LA}$-groupoid) is a Lie groupoid that is also a Lie algebroid; the unit manifold is also a Lie algebroid, and all the structure maps are Lie algebroid homomorphisms. A \emph{sub Lie algebroid-groupoid} is a subset that is both a Lie subalgebroid and a Lie subgroupoid.
\end{definition}
$T\G$ is the fundamental example of a Lie algebroid-groupoid. It is a Lie algebroid of the double groupoid $\Pair(\G)$. A sub Lie algebroid-groupoid is automatically a Lie algebroid-groupoid.
\begin{definition}
\label{wide}
A \emph{wide} subalgebroid \cite{mac3} is one that shares the same base manifold.
A \emph{full} subgroupoid \cite{mac3} is one that shares the same unit submanifold.
\end{definition}
In general, the base of a subalgebroid  is a submanifold. A wide subalgebroid is a subbundle in the usual sense, having a fiber over every point. This discrepancy between the terminologies for groupoids and algebroids is unfortunate. 

\begin{thm}
\label{restatement}
A real polarization of a Lie groupoid $\G$ is precisely a wide sub Lie algebroid-groupoid $\P\subset T\G$.
\end{thm}
\begin{proof}
An involutive real distribution is precisely a wide subalgebroid of the tangent bundle.

For a real polarization of a groupoid, the Hermiticity condition simplifies to 
\[
T\inv \P = \P .
\]
Multiplicativity and Hermiticity now mean that $\P\subset T\G$ is a subgroupoid. \end{proof}

The best behaved case of a real polarization $\P\subset T\G$ is one that integrates to a sub double groupoid of $\Pair(\G)$. (See \cite{b-m} for the definition of a double Lie groupoid.) This relates polarization to the following concept.
\begin{definition}
\label{congruence}
A \emph{smooth congruence} \cite{h-m2,mac3} of a Lie groupoid $\G$ is a closed, embedded sub double Lie groupoid of $\Pair(\G)$ that is full over $\G$.
\end{definition}
\begin{example}
If $G$ is a Lie group and $H\triangleleft G$ is a normal Lie subgroup, then the subset
\[
\{(g,g')\in\Pair(G) \mid g^{-1}g'\in H\}
\]
is a smooth congruence of $G$. Any smooth congruence of $G$ is of this form.
\end{example}
In this way, a congruence is a generalization of a normal subgroup. It is the most general structure for defining a quotient groupoid. As a full subgroupoid of $\Pair(\G)$, a congruence is an equivalence relation; the set of equivalence classes is the quotient groupoid.

The projection from a Lie groupoid to a quotient is a fibration. 
\begin{definition}
\label{fibration}
A \emph{fibration}  \cite{h-m2,mac3} of Lie groupoids is a smooth homomorphism $\Phi : \G\to \G'$ such that the base map $\Phi_0 : \G_0\to \G'_0$ and the map
\[
(\Phi,\s) : \G \longrightarrow \G' \mathbin{{}_\t\times_{\Phi_0}}\G_0
\]
are surjective submersions.
\end{definition}
Conversely, any fibration of Lie groupoids determines a smooth congruence as its ``kernel''. So, the notions of fibration and congruence are equivalent.

A real polarization $\P\subset T\G$ is in particular a foliation. If $\P$ integrates to a congruence, then the foliation is simple and the quotient groupoid is the leaf space. I will denote the quotient groupoid as $\G/\P$.

This suggests another definition, analogous to Definition \ref{admissible1} for complex polarizations of symplectic manifolds.
\begin{definition}
For a (complex) polarization $\P\subset T_\co\G$ of a groupoid,
two (generally singular) distributions $\D,\E\subseteq T\G$ are defined by 
$\D_\co :=\P\cap\bar\P$ and $\E_\co:=\P+\bar\P$. The polarization $\P$ is \emph{strongly admissible} if there exist groupoids and fibrations $\G \stackrel{p}\onto \G/\D \stackrel{q}\onto \G/\E$ such that $\D=\ker Tp$ and $\E=\ker T(q\circ p)$.
\end{definition}
In any case, if $\D$ has constant rank, then it is a real polarization of $\G$. If $\P$ is a strongly admissible polarization, then $\D$ and $\E$ are themselves strongly admissible real polarizations. However, if $\P$ is a symplectic groupoid polarization, then $\D$ and $\E$ will not be \emph{symplectic} groupoid polarizations, unless $\D=\E=\P$.

A strongly admissible real polarization is equivalent to a fibration whose fibers are connected. The kernel foliation $\P = \ker T\Phi$ of a fibration is always a strongly admissible real polarization. This is the easiest way to find groupoid polarizations.

As the example of a cotangent bundle showed, a polarized symplectic manifold can behave like its (unpolarized) quotient. This suggests that in the case of a strongly admissible real polarization, a polarized groupoid should act like its quotient. This will be the guiding principle in defining polarized convolution algebras. Not all polarizations are real and strongly admissible, but they'd like to be. A general polarized groupoid should be regarded as a virtual quotient.

\section{The Algebra}
\label{Algebra}
\subsection{Convolution}
There are two standard ways of constructing the convolution algebra of a Lie groupoid: using either a Haar system \cite{ren} or half-densities \cite{con1}. I favour the latter approach, because it is closer to the Hilbert space construction in geometric quantization. However, there is still an important discrepancy between half-forms and half-densities.

With this in mind, let's try to modify the definition of the convolution algebra by substituting half-forms for half-densities. This means using a bundle,
\beq
\label{half_form1}
\Omega^{1/2} := \sqrt{\Wedgem (T^{\t*}_\co\G\oplus T^{\s*}_\co\G)} .
\eeq
If the bundle $T^\t\G\oplus T^\s\G$ is orientable, then there exists a preferred ``positive'' choice of square root, and with this choice, half-forms are equivalent to half-densities. This orientability condition is frequently satisfied; it is true if the Lie algebroid $\Alg(\G)$ is an orientable bundle or if $\G$ is $\s$-simply connected.

If some other choice of square root is chosen, then it may still define a convolution algebra. However, it will not generally be possible to complete this to a \cs-algebra. The potential problem is exemplified by the $*$-algebra of matrices on an indefinite inner product space; that is not a \cs-algebra.

The convolution product of $a,b\in\Gamma_{\mathrm c}(\G,\Omega^{1/2})$ is defined by integrating $\pr_1^*a \pr_2^*b$ over the fibers of $\m$. For $\gamma\in\G$, the fiber $\m^{-1}(\gamma)$ is explicitly parametrized as
\beq
\label{m_fiber}
\m^{-1}(\gamma) = \{(\eta,\eta^{-1}\gamma) \mid \eta \in \t^{-1}(\t \gamma)\} .
\eeq
This gives the standard expression for the convolution product (see \cite{con1})
\beq
\label{convolution}
(a*b)(\gamma) := \int_{\eta\in\t^{-1}(\t\gamma)}  a(\eta) b(\eta^{-1}\gamma) .
\eeq
To make sense of this expression, note that \eqref{m_fiber} implies isomorphisms
\[
T^\m\G_2 \cong \pr_1^* T^\t\G \cong \pr_2^* T^\s\G .
\]
Using these, the coboundary bundle (see \eqref{L_coboundary}) is
\[
\partial^*\Omega^{1/2} \cong \Wedgem T^{\m*}_\co\G_2 \cong \pr_1^*\Wedgem T^{\t*}_\co\G
\]
the bundle of volume forms along the fibers of $\m$ in $\G_2$. With this isomorphism, the product $\pr_1^* a \pr_2^* b$ is a section of 
\[
\m^* \Omega^{1/2} \otimes  \pr_1^*\Wedgem T^{\t*}_\co\G .
\]
So, integrating this over $\t^{-1}(\t\gamma)$ in \eqref{convolution} does give another section of $\Omega^{1/2}$.

The other important property of $\Omega^{1/2}$ is the canonical isomorphism
\[
\inv^*\Omega^{1/2} \cong \overline{\Omega^{1/2}} .
\]
This defines the involution on the convolution algebra.

\subsection{Twisted Convolution}
\label{Twisted}
The concept of twisted \cs-algebras of a groupoid is fairly standard. Let $L$ be a line bundle over $\G$ and $\sigma$ a cocycle with coefficients in $L$. In \cite{ren}, Renault mainly discusses the case where $L$ is trivial, but he gives a more general definition: $\cs(\G,\sigma)$ is the quotient of $\cs(\G^\sigma)$ induced by the fundamental representation of $\T$, where $\G^\sigma$ is the extended groupoid in \eqref{extension}. In order to generalize this to the polarized case, it is better to work with sections of  $L\otimes \Omega^{1/2}$, and this is easy enough.

For $a,b\in \Gamma_{\mathrm c}(\G,L\otimes \Omega^{1/2})$, the twisted convolution product is defined by modifying \eqref{convolution} to
\[
(a*b)(\gamma) := \int_{\eta\in\t^{-1}(\t\gamma)} \sigma(\eta,\eta^{-1}\gamma)\, a(\eta)\, b(\eta^{-1}\gamma) .
\]
 
Although this is constructed with a cocycle $\sigma$, the algebra only depends upon the cohomology class of $[\sigma]\in \Tw(\G)$. 

\subsection{Polarized Convolution}
Let $\G$ be a Lie groupoid with a strongly admissible real polarization $\P\subset T\G$. Such a polarization is always the kernel foliation of a fibration 
\[
p: \G \onto \G/\P .
\] 
As I have suggested, the polarized groupoid should be treated as if it is the quotient groupoid. For this reason, I define the polarized convolution algebra to be
\[
\cs_\P(\G) := \cs(\G/\P)
\]
for any strongly admissible real polarization.
This case will be the guide to defining more general polarized convolution algebras. 

Let's reexpress this more directly in terms of $\G$. The elements of the polarized convolution algebra should be polarized sections of $p^*\Omega^{1/2}$. The pull-back of the target tangent bundle is
\[
p^* T^\t(\G/\P) \cong T^\t\G /(T^\t\G\cap \P) \mbox,
\]
so
\[
p^*\Omega^{1/2} \cong \Omega^{1/2}_\P := \sqrt{\Omega_\P} 
\]
is the ``positive'' square root of 
\[
\Omega_\P := \Wedgem\left(T^\t\G/[T^\t\G \cap \P]\oplus T^\s\G / [T^\s\G \cap \P]\right)^*_\co.
\]
The key property of this bundle is that 
\[
\partial^*\Omega^{1/2}_\P \cong \pr_1^* \Wedgem(T^\t\G/[T^\t\G \cap \P])^* .
\]
This is (the pull-back of) the bundle of volume forms along the $\t$-fibers but transverse to $\P$. If $a,b\in \Gamma(\G,\Omega^{1/2}_\P)$ are two global polarized sections (i.e., the pull-backs of sections over $\G/\P$) and are compactly supported modulo $\P$ then the polarized convolution is defined by modifying the integral \eqref{convolution}. Rather than integrating over the fiber $\t^{-1}[\t(\gamma)]$, we can note that the integrand is constant along $\P$-fibers, so it can be integrated over the quotient space. Alternately, we can integrate over a transversal submanifold-with-boundary of $\t^{-1}[\t(\gamma)]$ that intersects each $\P$-leaf there exactly once.

For the moment, this is just an awkward reformulation of convolution on the quotient groupoid $\G/\P$. However, it is justified by its generalizations.

For a complex polarization, $\Omega^{1/2}_\P$ generalizes trivially to
\beq
\label{half_form}
\begin{split}
\Omega^{1/2}_\P &:= \sqrt{\Omega_\P} \\
\Omega_\P &:= \Wedgem\left(T^\t_\co\G/[T^\t_\co\G \cap \P]\oplus T^\s_\co\G / [T^\s_\co\G \cap \P]\right)^*
\end{split}
\eeq
provided that $T^\t\G\cap\P$ has constant rank. This gives a general definition for $\Omega^{1/2}_\P$, if $T^\t\G\cap\P$ has constant rank and there exists a square root bundle satisfying 
\[
\inv^*\Omega^{1/2}_\P \cong \overline{\Omega^{1/2}_\P} .
\]
This introduces a possible ambiguity because of the choice of square root; any other choice is obtained by tensoring $\Omega^{1/2}_\P$ with a real line bundle which is isomorphic to its pull-back by $\inv$. This choice is parametrized by the $\inv$-invariant cohomology, $H^1(\G;\Z_2)^{\inv}$. However, it is plausible that only one choice will lead to \cs-algebras.

For a strongly admissible real polarization, a section of $\Omega^{1/2}_\P$ is polarized if it is the pull-back by $p$ of a section over $\G/\P$. This can also be expressed in terms of a flat $\P$-connection on $\Omega^{1/2}_\P$. Showing that this connection exists more generally requires a couple of steps.

\begin{definition}
For any distribution $\P\subset T_\co\G$, define the (generally singular) distribution $\P_0\subset T_\co M$ by
\[
\P_{0\, x} := T\t(\P_x)
\]
for every $x\in M$.
\end{definition}
\begin{lem}
\label{P0.lem}
If $\P$ is a multiplicative distribution, then for any $\gamma\in\G$, 
\[
T\t(\P_\gamma)= \P_{0\;\t(\gamma)} .
\]
If $\P$ is also involutive, then $\P_0$ is involutive.
\end{lem}
\begin{proof}
The condition of multiplicativity at $(\gamma,\gamma^{-1})\in\G_2$ states that,
\[
\P_{\t(\gamma)} = T\m\left[\P_{2\, (\gamma,\gamma^{-1})}\right]  .
\]
Applying $T\t$, and using $T\t[X\cdot Y]=T\t(X)$, this gives
\[
\begin{split}
\P_{0\:\t(\gamma)} &= \left\{T\t(X)\bigm| X\in\P_\gamma, \exists Y\in\P_{\gamma^{-1}}: T\s(X)=T\t(Y)\right\}\\
 &\subseteq T\t(\P_\gamma) .
\end{split}
\]
Conversely, multiplicativity at $(\t(\gamma),\gamma)$ says,
\[
\P_\gamma = T\m\left[\P_{2\,(\t\gamma,\gamma)}\right] 
\]
and implies
\[
\begin{split}
T\t(\P_\gamma) &= \left\{T\t(X)\bigm| X\in\P_{\t\gamma}, \exists Y\in\P_\gamma: T\s(X)=T\t(Y)\right\}\\
&\subseteq T\t(\P_{\t\gamma}) = \P_{0\;\t(\gamma)} .
\end{split}
\]

Now suppose that $\P$ is involutive and recall that a vector field is "projectable" by a smooth map if the push-forward is well defined.
Any section of $\P_0$ can be lifted to a $\t$-projectable section of $\P$. If $X,Y\in\Gamma(\G,\P)$ are $\t$-projectable, then $[X,Y] \in \Gamma(\G,\P)$ and so
\[
[T\t(X),T\t(Y)] = T\t([X,Y])
\]
is a section of $\P_0$. Therefore $\P_0$ is involutive.
\end{proof}

\begin{thm}
If $\P$ is a polarization of a groupoid $\G$, then the $\P$-Bott connection (Def.~\ref{Bott}) induces a natural flat $\P$-connection on $\Omega_\P$ and $\Omega^{1/2}_\P$, provided that these bundles are defined.
\end{thm}
\begin{proof}
The inverse image of $\P_0$ is $T\t^{-1}(\P_0) = \P+ T^\t_\co\G$ and the $\C^\infty(\G)$-module of sections is generated by the $\t$-projectable ones. So, let $X,Y\in\Gamma(\G,\P+ T^\t_\co\G)$ be $\t$-projectable. Then $T\t([X,Y])$ is a section of $\P_0$, so $[X,Y]$ is a section of $\P+ T^\t_\co\G = T\t^{-1}(\P_0)$, therefore $\P+ T^\t_\co\G$ is involutive.

This implies that the Bott connection for $\P$ preserves $[\P+ T^\t_\co\G ]^\perp \subset \P^\perp$. So, this induces a natural flat $\P$-connection on
\[
\P^\perp / [\P+ T^\t_\co\G ]^\perp 
\cong \left[(\P+ T^\t_\co\G )/\P\right]^*
\cong \left[T^\t_\co\G/(\P\cap T^\t_\co\G)\right]^* .
\]
The same is true with $T^\s\G$ in place of $T^\t\G$, so this gives the desired $\P$-connection on $\Omega_\P$. Finally, the square root $\Omega^{1/2}_\P$ inherits a connection by a simple application of the Leibniz rule.
\end{proof}

With this connection, we can speak of polarized sections of $\Omega^{1/2}_\P$. If $U\subset \G$ is an open subset, then $a \in \Gamma(U,\Omega^{1/2})$ is polarized if $\nabla_X a =0$ for any $X\in\P$ over $U$. 

Some of the examples later in this paper suggest that it is useful to consider polarizations for which $T^\t_\co\G\cap\P$ does not have constant rank. In that case, $T^\t_\co\G\cap\P$  is not a bundle, so $\Omega^{1/2}_\P$ is not defined. However, I expect that the sheaf of local polarized sections of $\Omega^{1/2}_\P$ may be meaningful even when the bundle is not.

Now suppose that $\P\subset T_\co\G$ is a strongly admissible polarization and $\Omega^{1/2}_\P$ does exist.
Given two polarized sections $a,b\in\Gamma(\G,L\otimes\Omega^{1/2}_\P)$, the convolution $(a*b)(\gamma)$ should be some sort of integral of 
\beq
\label{outer_product}
 a(\eta)\, b(\eta^{-1}\gamma)
\eeq
with $\t(\eta)=\t(\gamma)$. 
The polarization of $a$ and $b$ implies that \eqref{outer_product} is covariantly constant along the foliation $\D$ restricted to $\t^{-1}[\t(\gamma)]$, so we should integrate over the quotient by $\D$.
The problem is that for such an integration, we need an isomorphism 
\beq
\label{product_isomorphism}
\partial^*\Omega^{1/2}_\P  \isom  \pr_1^*\Wedgem(T^\t\G/[T^\t\G\cap \D])^*_\co .
\eeq
Instead, there is a natural isomorphism
\[
\partial^*\Omega_\P \cong \pr_1^*\Wedgem(T^\t\G/[T^\t\G\cap \D])^*_\co \otimes \pr_1^*\Wedgem(T^\t\G/[T^\t\G\cap \E])^*_\co
\]
To ``correct'' this, we need an isomorphism, 
\beq
\label{correction}
\Wedgem(T^\t\G\cap \E) \isom \Wedgem(T^\t\G\cap \D) .
\eeq
This requires some additional structure. I will show in the next section that a symplectic groupoid has such a structure.

I will not go any further in trying to define polarized convolution for complex polarizations without a twist, because of two problems. One is the need for an isomorphism \eqref{correction}. 
The other problem is that suitable global polarized sections may not exist, as illustrated by this example.
\begin{example}
Consider the pair groupoid $\Pair(\T^2)$ of a complex torus with the complex structure reversed on the second factor. Let $\P$ be the antiholomorphic tangent bundle. The bundle $\Omega^{1/2}_\P$ is trivial in this case, so global polarized sections of $\Omega^{1/2}_\P$ are just holomorphic functions on $\T^4$, but only constant functions are holomorphic. Global polarized sections of $\Omega^{1/2}_\P$ are not suitable for defining a convolution algebra.

On the other hand, if this is twisted by a line bundle with positive curvature, then many suitable holomorphic sections will exist.
\end{example}
This leads to the final level of generalization.

\subsection{Twisted Polarized Convolution}
I will concentrate on the most relevant case: that of a symplectic groupoid. However, twisted polarized convolution should also make sense in some other cases, including at least real polarizations and some polarizations of Poisson groupoids.

\begin{thm}
If $\P$ is a polarization of a symplectic groupoid $\Sigma$ and $\rk T^\t_\co\Sigma\cap \P$ is constant, then for some $k\in\mathbb{N}$, contraction with $\frac{\omega^k}{k!}$ gives an isomorphism,
\beq
\label{correction2}
\Wedgem(T^\t\Sigma\cap \E) \isom\Wedgem(T^\t\Sigma\cap \D) .
\eeq
\end{thm}
\begin{proof}
The tangent bundle restricted to the unit submanifold $\unit: M\into\Sigma$ is a subgroupoid $\unit^* T\Sigma \subset T\Sigma$. 
It is easy to explicitly construct the groupoid structure on $\unit^*T\Sigma$ in terms of the anchor map. In particular, the groupoid multiplication on $\unit^*(T^\t\Sigma\cap T^\s\Sigma)$ is just addition, so for any $x\in M$ and $X\in T^\t_x\Sigma\cap T^\s_x\Sigma$, we have $X^{-1}=-X$.

Now, for any element $\gamma\in\Sigma$ and any vector $X\in T^\t_{\co\;\gamma}\Sigma\cap T^\s_{\co\;\gamma}\Sigma\cap\P_\gamma$, Hermiticity implies $X^{-1}\in\bar\P$. Multiplicativity implies that $0_\gamma\cdot X^{-1}\cdot 0_\gamma\in \P$, where $0_\gamma\in T_\gamma\Sigma$ is the $0$ vector over $\gamma$ and a dot denotes the tangent groupoid multiplication, $T\m$. However, 
\[
0_\gamma\cdot X^{-1}\cdot 0_\gamma = (X\cdot 0_{\gamma^{-1}})^{-1} \cdot 0_\gamma = (-X\cdot0_{\gamma^{-1}})\cdot 0_\gamma = -X .
\]
So, $T^\t_\co\Sigma\cap T^\s_\co\Sigma\cap \P = T^\t_\co\Sigma\cap T^\s_\co\Sigma\cap \bar\P$. This means that $T^\t\Sigma\cap T^\s\Sigma\cap \E = T^\t\Sigma\cap T^\s\Sigma\cap \D$.

Recall that the symplectic orthogonal bundle to $T^\t\Sigma$ is $\#(T^\t\Sigma)^\perp = T^\s\Sigma$, and the symplectic orthogonal to $\E$ is $\D$. If we restrict the symplectic form $\omega$ to $T^\t\Sigma\cap \E$, then the kernel is 
\[
T^\t\Sigma\cap \E \cap \#(T^\t\Sigma\cap \E)^\perp = T^\t\Sigma \cap (T^\s\Sigma\cap\E + \D) = T^\t\Sigma\cap \D .
\]
This shows that $\omega$ gives a \emph{nondegenerate} form on $(T^\t\Sigma\cap\E)/(T^\t\Sigma\cap\D)$. The dimension of this bundle must be even, say $2k$, and $\frac{\omega^k}{k!}$ gives the desired isomorphism.
\end{proof}

Let $\P$ be a polarization of a symplectic groupoid $\Sigma$ and $(\sigma,L,\nabla)$ a prequantization. The connection on $L$ combines with the partial connection on $\Omega^{1/2}_\P$ to give a flat $\P$-connection on $L\otimes\Omega^{1/2}_\P$. So, local polarized sections will exist.

Suppose that the polarization $\P$ is strongly admissible.
If $a,b\in \Gamma(\Sigma,L \otimes \Omega^{1/2}_\P)$ are polarized sections, then the compatibility of $\sigma$ with the connection implies that 
\beq
\label{outer_product2}
\sigma(\eta,\eta^{-1}\gamma)\, a(\eta)\, b(\eta^{-1}\gamma)
\eeq
is covariantly $\D$-constant as a function of $\eta\in \t^{-1}[\t(\gamma)]$. To define the twisted convolution $a*b$, we should integrate \eqref{outer_product2} over the quotient. If the square root $\Omega^{1/2}_\P$ can be chosen appropriately, then the isomorphism \eqref{correction2} induces an isomorphism
\[
\partial^*\Omega^{1/2}_\P \cong \pr_1^*\Wedgem(T^\t\Sigma/[T^\t\Sigma\cap \D])^*_\co
\]
which makes this integration meaningful.

There are two subtleties in defining a convolution algebra from polarized sections of $L\otimes \Omega^{1/2}_\P$. The first is the problem of fall-off conditions. Without a polarization, a convolution algebra is usually defined with compactly supported sections, although a somewhat weaker fall-off condition would produce the same \cs-algebra in the end. For a strongly admissible polarization, compactly supported polarized sections will only exist if the leaves of the foliation $\E$ are compact. It is a delicate matter to formulate a general fall-off condition that is strong enough, but not too strong. I do not solve this problem here.

The second issue is that not all polarizations are strongly admissible, and even for strongly admissible polarizations the space of polarized sections may be ``too small''.
The idea of ``cohomological wave functions'' suggests that we should not just use the global polarized sections, but rather the convolution algebra should be the total cohomology of the sheaf of (local) polarized sections of $L\otimes\Omega^{1/2}_\P$. Unfortunately, the problem of defining the inner product of cohomological wave functions translates here into the problem of defining convolution on cohomology spaces. I do not attempt to solve this here either.

If $\P$ is strongly admissible, $\Omega^{1/2}_\P$ is a bundle, and the higher degree cohomology of polarized sections of $L\otimes \Omega^{1/2}_\P$ vanishes, \emph{then} the convolution algebra should consists of global polarized sections.

\subsection{Completion}
\label{Completion}
Once a convolution algebra has been constructed, the final step is to complete this to a \cs-algebra, $\cs_\P(\Sigma,\sigma)$.

Since this algebra is a generalization of the convolution algebra of a groupoid, we can expect that there will be more than one natural way of completing it. 
Formally, the easiest to define will be the maximal completion; this is defined using (essentially) all possible representations.
The reduced \cs-algebra should be defined using a ``regular'' representation.

In the correspondence between Poisson manifolds and \cs-algebras, symplectic manifolds correspond to Hilbert spaces (or sometimes, Hilbert modules) and Poisson maps correspond (contravariantly) to $*$-homomorphisms. Based on this, symplectic realizations should correspond to representations.

Indeed, the regular representation (which I have not defined) should correspond to $\t:\Sigma\to M$, which is itself a symplectic realization of $M$.
There may be a general procedure for quantizing symplectic realizations to representations. If so, then this gives a preferred class of ``geometrical'' representations. Completion using these representations may give a third way --- a geometrical \cs-algebra completion.

By convention, $\cs(\G)$ denotes the maximal \cs-algebra; I am not using the notation $\cs_\P(\Sigma,\sigma)$ in this sense, but rather out of agnosticism. I do not know which completion will fit a reasonable definition of quantization.

Some examples (Sec.~\ref{Linear}) suggest that both the maximal and reduced \cs-algebras are quantizations. One example (Sec.~\ref{Multiply}) suggests that the reduced \cs-algebra is the natural quantization. The results in \cite{h-l} suggest that the maximal \cs-algebra is the best behaved with respect to symplectic reduction.

\subsection{Real Polarizations and Bohr-Sommerfeld Conditions}
Suppose that $\P$ is a strongly admissible real polarization of a symplectic groupoid $\Sigma$. Let $p:\Sigma\onto\Sigma/\P$ be the quotient fibration. By design, $\cs_\P(\Sigma)\cong \cs(\Sigma/\P)$, but with a twist, things become more subtle.

Let $(\sigma,L,\nabla)$ be a prequantization of $\Sigma$. If the leaves of $\P$ are simply connected, then the connection canonically trivializes $L$ along these leaves. This identifies $L$ as the pull back $L\cong p^* L_0$ of some bundle over $\Sigma/\P$. The algebra is thus,
\[
\cs_\P(\Sigma,\sigma) \cong \cs(\Sigma/\P,\sigma_0)
\]
where $\sigma_0$ is a cocycle (the ``reduced cocycle'') with coefficients in $L_0$.  Note that $L_0$ does not inherit a connection.

To compute the algebra up to isomorphism, we don't need the specific cocycle $\sigma_0$, just its cohomology class $[\sigma_0] \in \Tw(\Sigma/\P)$. 
The twist group $\Tw$ is a contravariant functor, so it maps $\Tw(p):\Tw(\Sigma/\P)\to\Tw(\Sigma)$. The reduced twist is related to the prequantization by $\Tw(p)(\sigma_0) \simeq \sigma$.

In most of the examples, $L$ is actually trivializable. A trivialization identifies $L$ with $\co\times\Sigma$. The connection is then given by a $1$-form $\theta\in\Omega^1(\Sigma)$ as $\nabla= d +i\theta$. The properties of a prequantization imply that $d\theta=-\omega$ and $\unit^*\theta$ is exact. We can always apply a gauge transformation to get another trivialization for which $\unit^*\theta=0$. 
\begin{definition}
A \emph{symplectic potential} is a $1$-form $\theta\in\Omega^1(\Sigma)$ such that $d\theta=-\omega$ and $\unit^*\theta=0$. It is \emph{adapted} if it is conormal to the polarization, i.e., $\theta\in\Gamma(\Sigma,\P^\perp)$.
\end{definition}

If there exists an adapted symplectic potential, then it is easy to compute the reduced cocycle.
\begin{lem}
Let $\P$ be a strongly admissible real symplectic groupoid polarization with simply connected leaves, and $p:\Sigma\to\Sigma/\P$ the projection to the quotient. If $\theta$ is an adapted symplectic potential, then $L_0$ is trivial and the reduced cocycle $\sigma_0\in \C^\infty([\Sigma/\P]_2,\T)$ is determined (up to a locally constant phase) by
\[
p^*(\sigma_0^{-1}d\sigma_0) = i\,\partial^*\theta.
\]
\end{lem}
\begin{proof}
The simplicial coboundary is a linear combination of pullbacks, so it commutes with the exterior derivative. Recall that the simplicial coboundary $\partial^*\theta$ measures the failure of $\theta$ to be multiplicative. However, the symplectic form is multiplicative, so 
\[
d\partial^*\theta=\partial^* d\theta = -\partial^*\omega = 0 \mbox,
\]
and thus $\partial^*\theta$ is closed. 

Recall that the cocycle $\sigma$ is a covariantly constant section of the coboundary bundle $\partial^*L^*$ over $\Sigma_2$. Existence of a symplectic potential implies that this bundle is trivial and $-\partial^*\theta$ is its connection $1$-form. So,
\[
0 = \nabla \sigma = d\sigma - i (\partial^*\theta)\sigma
\]
and
\[
\sigma^{-1} d\sigma = i\,\partial^*\theta .
\]

By assumption, $\theta$ is conormal to $\P$, so $\partial^*\theta$ is conormal to $\P_2\subset T\Sigma_2$. This means that $\sigma$ is $\P_2$-constant, so it is the pullback of some $\sigma_0\in\C^\infty([\Sigma/\P]_2,\T)$.

By definition, the sections of the reduced line bundle $L_0$ are identified with the covariantly $\P$-constant sections of $L$. The assumption that $\theta$ is adapted means that ``covariantly $\P$-constant'' is just $\P$-constant in this trivialization. So, $L_0$ is trivial. With these trivializations, the reduced cocycle is the function $\sigma_0$.
\end{proof}
In practice, there usually exists a real function $\phi\in\C^\infty([\Sigma/\P]_2)$ such that 
\[
d\,p^*\phi = \partial^*\theta
\]
and so $\sigma_0=e^{i\phi}$.

\begin{cor}
\label{no_twist}
If $\Sigma$ is $\s$-connected and admits a multiplicative adapted symplectic potential, then the reduced cocycle is trivial.
\end{cor}
\begin{proof}
Multiplicativity means that $\partial^*\theta=0$, so $\sigma_0$ is locally constant
\end{proof}

If the leaves of $\P$ are not simply connected, then $L$ will typically have holonomy around them. In this case, there is an open set over which smooth polarized sections of $L\otimes\Omega^{1/2}_{\P}$ must vanish, and it is clearly insufficient to construct a convolution algebra from smooth globally polarized sections.

This is the same situation encountered in the construction of a Hilbert space. In that case, the problem was tentatively solved with the idea of distributional wave functions and Bohr-Sommerfeld quantization. In \cite{wei3}, Weinstein derived the noncommutative $2$-torus by applying Bohr-Sommerfeld quantization to a symplectic groupoid. 
\begin{definition}
The \emph{Bohr-Sommerfeld subgroupoid} $\Sigma_{\BS}\subset\Sigma$ is the set of points through which $L\otimes\Omega^{1/2}_\P$ has trivial holonomy. The \emph{reduced groupoid} is the quotient $\Sigma_\BS/\P$.
\end{definition}
This terminology is deliberate. The construction of the ``reduced groupoid'' is very similar to symplectic reduction.

Tentatively, the algebra should be 
\[
\cs_\P(\Sigma,\sigma) \cong \cs(\Sigma_\BS/\P,\sigma_0)
\]
where $(\sigma_0,L_0)$ is an inherited twist.

Let $\iota$ be the inclusion of $\Sigma_\BS$ and $p$ the quotient map:
\[
\Sigma \stackrel{\iota}{\hookleftarrow} \Sigma_\BS \stackrel{p}{\onto} \Sigma_\BS/\P .
\]
The reduced twist is related to the prequantization by $\Tw(\iota)(\sigma)\simeq \Tw(p)(\sigma_0)$.

More generally, if $\P$ is a (complex) strongly admissible polarization, then the Bohr-Sommerfeld conditions come from holonomy around the leaves of $\D$. There is again a Bohr-Sommerfeld subgroupoid $\Sigma_\BS$ and a reduced groupoid $\Sigma_\BS/\D$. The twisted polarized \cs-algebra of $\Sigma$ should be a twisted polarized \cs-algebra of the reduced groupoid, but now by a totally complex polarization. In general, the reduced groupoid is a Poisson groupoid.

\section{Examples}
\label{Examples}
\subsection{Symplectic}
This first example was already used to motivate the definitions. If $M$ is a symplectic manifold with a polarization $\F$, then $\P:=\F\times\bar\F$ is a symplectic groupoid polarization of $\Pair(M)$.

For the pair groupoid $\Sigma = \Pair(M)$, the set of composable pairs can be identified as $\Sigma_2\cong M\times M\times M$. The multiplication map $\m:\Sigma_2\to\Sigma$ simply forgets the middle factor.  Now $\P_2 = (\P\times\P)\cap T_\co\Sigma_2 = \F\times(\F\cap\bar\F)\times \bar\F$, so the multiplicativity is apparent.

I shall call this type of induced polarization \emph{exact}, because it can be thought of as the ``coboundary'' of the symplectic polarization $\F$. Not all polarizations of symplectic pair groupoids are exact; I give such an example in the next section.

\begin{example}
Suppose that $M=T^*N$ for some smooth manifold $N$. The natural polarization of $T^*N$ is the ``vertical'' polarization, whose leaves are the cotangent fibers. The pair groupoid is $\Pair(T^*N) \cong T^*(\Pair N)$. The vertical polarization is just the kernel foliation of the projection down to $\Pair(N)$, so that is the reduced groupoid.

We can take the Liouville $1$-form as a symplectic potential. This is adapted to the vertical polarization, but it is also multiplicative, therefore the quantization is the untwisted groupoid \cs-algebra 
\[
\cs(\Pair N) \cong \K[L^2(N)]
\mbox,
\]
as expected.
\end{example}

A real polarization $\F$ of $M$ induces the exact polarization $\P=\F\times\F$, the set of pairs of vectors from $\F$. Two such pairs are composable if the second of the first pair is the first of the second pair, so $\P_2 = \F\times \F\times \F$.

The case of a purely complex polarization lies at the opposite extreme.
\begin{example}
Let $M$ be a \Kahler\ manifold with $\F$ the antiholomorphic tangent bundle. This induces the exact polarization $\P:=\F\times\bar\F$. Two pairs of vectors are only composable if the second of the first and the first of the second pair both vanish. That is, $\P_2 = \F\times 0 \times \bar \F$.

This polarization is equivalent to a \Kahler\ structure on the groupoid itself. In this case, the convolution product is (essentially) ordinary convolution of sections which happen to be polarized (holomorphic).
\end{example}

More generally, suppose that $M$ is a symplectic manifold with polarization $\F$ and prequantization $L$. Let $\Sigma:=\Pair(M)$, $\P:=\F\times\bar\F$, and $\sigma$ be the prequantization with bundle $L\boxtimes\bar L$. To see that this
gives $\cs_\P(\Sigma,\sigma)= \K[L^2_\F(M,L)]$, first note that $T^\t_\co\Sigma \equiv \ker \t = 0 \times \TcM = \s^* \TcM$. So 
\[
T^\t_\co \Sigma / (T^\t_\co\Sigma \cap \P) = \s^* (\TcM/\bar\F) = \s^* (\bar\F^\perp)^*
\]
and likewise $T^\s_\co \Sigma / (T^\s_\co\Sigma \cap \P) = \t^*(\F^\perp)^*$. Referring to eq.~\eqref{half_form}, this shows that 
\[
\Omega_\P = \Wedgem \F^\perp \boxtimes \Wedgem\bar\F^\perp
\]
and the square root can be chosen so that the half-form bundle is
\[
\Omega_\P^{1/2} = \sqrt{\Wedgem \F^\perp} \boxtimes \sqrt{\Wedgem\bar\F^\perp} .
\]
Which is exactly what we wanted. 
The algebra $\cs_\P(\Sigma,\sigma)$ is by definition constructed from polarized sections of $L\otimes \sqrt{\Wedgem \F^\perp} \boxtimes \bar L\otimes \sqrt{\Wedgem \bar\F^\perp}$. This space (or cohomology) of polarized sections over $\Sigma=M\times M$ is just a tensor product of the space (or cohomology) of polarized sections of $L\otimes \sqrt{\Wedgem \F^\perp}$ over $M$ with its complex conjugate. 
To the extent that either algebra is defined, we do have $\cs_\P(\Sigma,\sigma)= \K[L^2_\F(M,L)]$.

\subsection{Constant}
\label{Constant}
Suppose that our Poisson manifold is a vector space $V$ with a constant Poisson bivector $\pi\in\Wedge^2 V$. The symplectic integration is $\Sigma= V\oplus V^*$ as a vector space and  $\Sigma =T^*V$, symplectically. All the structure maps are linear.

Let $x^i$ be coordinates on $V$ and $y_i$ coordinates on $V^*$. The source and target maps are
\[
\s(x^i,y_i) = x^i - \tfrac12 \pi^{ji}y_j 
\]
\[
\t(x^i,y_i) = x^i + \tfrac12 \pi^{ji}y_j .
\]
The set of composable pairs can be identified with $\Sigma_2\cong V\oplus V^*\oplus V^*$; the relevant maps are
\begin{gather*}
\pr_1(x^i,y_i,y'_i) = (x^i + \tfrac12 \pi^{ji}y'_j,y_i) \\
\pr_2(x^i,y_i,y'_i) = (x^i - \tfrac12 \pi^{ji}y_j,y'_i) \\
\m(x^i,y_i,y'_i) = (x^i,y_i+y'_i) .
\end{gather*}

The projection from the groupoid $\Sigma\cong V\oplus V^*$ to the additive group $V^*$ is a fibration of groupoids. The fibers are Lagrangian, so this is a polarization of the symplectic groupoid. The simplest choice of adapted symplectic potential is
\[
\theta = -x^idy_i .
\]
We can compute directly
\begin{align*}
\partial^* \theta &= \left(x^i+\tfrac12\pi^{ji}y'_j\right)dy_i - x^i\left(dy_i+dy'_i\right) + \left(x^i-\tfrac12\pi^{ji}y_j\right)dy'_i \\
&= -\tfrac12\pi^{ij}y'_jdy_i - \tfrac12\pi^{ij}y_idy'_j \\
&= d\left(-\tfrac12\pi^{ij}y_iy'_j\right) .
\end{align*}
This gives a \emph{group} cocycle $\sigma_0:V^*\times V^*\to\T$,  $\sigma_0(y,y')= e^{-\tfrac i2\pi(y,y')}$. So, my quantization recipe applied to  $(V,\pi)$ gives the twisted group algebra, $\cs(V^*,\sigma_0)$.

This is the usual Moyal quantization of a Poisson vector space, see \cite{rie3}. 
When $\pi=0$, the Fourier transform gives the canonical isomorphism $\cs(V^*)\cong \C_0(V)$. Because both $\cs(V^*)$ and $\cs(V^*,\sigma_0)$ are completions of the same vector space (compactly supported smooth densities) we can identify dense subalgebras as vector spaces. With these identifications, this quantization gives a deformed product on $\C^\infty_0(V)$. The product is given by an integral kernel which is the Fourier transform of $\sigma_0$.

\begin{example}
If $V$ is a symplectic vector space, then $\Sigma\cong \Pair(V)$ is the same symplectic groupoid, but this ``horizontal'' polarization does not come from a symplectic polarization of $V$. Nevertheless, the algebras given by these polarizations are isomorphic $\K[L^2(V)] \cong \cs(V^*,\sigma_0)$ because $\pi$ is nondegenerate.
\end{example}

\subsection{Linear}
\label{Linear}
Suppose that $M$ is a vector space with a Poisson bivector that is linear. Then $M$ can be identified with the linear dual $\mathfrak g^*$ of a Lie algebra, with the Lie-Poisson structure.

Let $G$ be a connected Lie group integrating $\mathfrak g$. We can give the cotangent space $T^*G$ a symplectic groupoid structure integrating $\mathfrak g^*$; see \cite{c-d-w,mac3}. First, let the unit be the map $\unit:\mathfrak g^*\to T^*_eG\subset T^*G$ identifying $\mathfrak g^*$ with the cotangent fiber over the group unit. Then, define the source and target maps as the duals of the maps which identify $\mathfrak g$ with left and right-invariant vector fields. Finally, define the multiplication to be the unique one compatible with these maps such that the product of an element of $T^*_gG$ with one of $T^*_hG$ lies in $T^*_{gh}G$.

The bundle projection $p:T^*G\onto G$ is a fibration of groupoids, therefore $\ker Tp$ is a real polarization of the groupoid $T^*G$. The fibers are Lagrangian, therefore this is a symplectic groupoid polarization. The Liouville $1$-form is a multiplicative adapted symplectic potential, therefore the reduced twist will be trivial. The leaves of this polarization are the cotangent fibers, which are contractible, so there are no Bohr-Sommerfeld conditions. The reduced groupoid is therefore the quotient groupoid, which is the group $G$. The quantization recipe applied to $\mathfrak g^*$ gives either $\cs_{\ker Tp}(T^*G,\sigma)=\cs(G)$ or $\csr(G)$, depending upon the choice of completion (see Sec.~\ref{Completion}).

This is a standard example of quantization. It was first proposed by Rieffel \cite{rie4}, who proved that with the reduced \cs-algebra $\csr(G)$, there is a ``deformation quantization by partial embeddings''. In the case that $G$ is nilpotent, he proved that $\cs(G)$ gives a strict deformation quantization. 

The restriction to nilpotent groups is related to a conjecture of Weinstein \cite{wei3,wei4} that a ``good quantization'' of a Poisson manifold exists if and only if it is of ``exponential type'', meaning that there exists a symplectomorphism $\Sig(M)\cong T^*M$ that intertwines the groupoid inverse with negation.
The dual of a nilpotent Lie algebra is of exponential type. The dual of $\mathfrak{su}(2)$ is not, because $T^*S^3$ is not homeomorphic to $\R^6$, let alone symplectomorphic.

These notions of ``strict'' and ``nice'' are based on actually identifying dense subalgebras of the quantum and classical algebras as vector spaces. In the terms of Section~\ref{Quantization}, this means requiring injective quantization maps. 
In my view, this is quite unnecessary and overly restrictive. $\cs[\mathrm{SU}(2)]$ is a very nice quantization indeed!

If $A$ is a Lie algebroid over a manifold $N$, then (the total space of) the dual bundle $A^*$ is a Poisson manifold \cite{c-d-w,cou}. This example includes both $\mathfrak g^*$ and $T^*N$ as special cases. A Poisson manifold is of this form if it is a vector bundle and the space of smooth fiber-linear functions is closed under the Poisson bracket.

If $\G$ is a groupoid integrating $A$, then the cotangent space $T^*\G$ has a natural structure as a symplectic groupoid over $A^*$, see \cite{c-d-w,mac3}. 

Let $p:T^*\G\to\G$ be the bundle projection. This is a fibration of groupoids, therefore $\ker Tp$ is a groupoid polarization. This is also a Lagrangian foliation, therefore $\ker Tp$ is a symplectic groupoid polarization. Once again, $T^*\G$ is a cotangent space and the Liouville form is a multiplicative adapted symplectic potential; therefore the twist can be eliminated by Corollary~\ref{no_twist}. The quantization of $A^*$ is thus 
\[
\cs_{\ker Tp}(T^*\G,L) \cong \cs(\G) .
\] 

Landsman \cite{lan2} proposed this as an example of quantization. There is a parallel result for formal quantization in \cite{n-w-x}. Landsman and Ramazan \cite{l-r} proved that \emph{both} $\cs(\G)$ and $\csr(\G)$ give strict deformation quantizations of $A^*$, although their definition of ``strict'' differs from Rieffel's.
This suggests that in my general construction \emph{all} completions may be legitimate quantizations.

\subsection{Multiply Connected}
\label{Multiply}
Let $M$ be a connected symplectic manifold with universal covering $\Mt$ and fundamental group $\Gamma:=\pi_1(M)$. Given that $M$ is symplectic, the most obvious symplectic integration is $\Pair(M)$. However, this integration does not generalize to other Poisson manifolds. For an arbitrary integrable Poisson manifold, the canonical symplectic integration is the (unique) $\s$-simply connected symplectic groupoid. In this case, that is the fundamental groupoid 
\[
\Sig(M) = \PI(M) \cong \Pair(\Mt)/\Gamma \mbox,
\]
where $\Gamma$ acts diagonally. This is the $\s$-simply connected covering of $\Pair(M)$. 

As above, any polarization $\F$ of the symplectic manifold $M$ determines a polarization of $\Pair(M)$, but because $\PI(M)$ is a covering of $\Pair(M)$, $\F$ also gives a polarization of $\PI(M)$.

$\PI(M)$ may be prequantizable when $M$ is not. Any prequantization of $\PI(M)$ can be constructed from a prequantization of $\Mt$ with a projective action of $\Gamma$. This exists if and only if the integral of $\omega$ over any $S^2$ is a multiple of $2\pi$, as expected from Theorem~\ref{c_prequantizability}.

This quantization is heuristically equivalent to quantizing $\Mt$, and then quotienting by $\Gamma$. I studied such a construction in \cite{haw5} for $M$ compact and \Kahler. A Toeplitz quantization map can be constructed for $\Mt$ but applied to $\C(M)$. The \cs-algebra generated by the image is a quantization of $M$. I showed that this \cs-algebra is isomorphic to the algebra of compact operators on a Hilbert $\csr(\Gamma,\sigma_\Gamma)$-module, where the group cocycle $\sigma_\Gamma$ is derived from the symplectic form. 
I expect that this algebra is isomorphic to $\cs_{\P,\mathrm r}[\PI(M),L]$.
The best way to prove this would be to show that $\cs_{\P,\mathrm r}[\PI(M),L]$ is the algebra of compact operators on a Hilbert $\csr(\Gamma,\sigma_\Gamma)$-module and prove that these Hilbert modules are isomorphic. It is plausible that such a Hilbert module can be constructed by (in some sense) quantizing the Morita equivalence between $\PI(M)$ and $\Gamma = \pi_1(M)$.

In this example, it is specifically the reduced \cs-algebra, $\csr(\Gamma,\sigma_\Gamma)$ that emerges. This suggests that the reduced \cs-algebra may be the most relevant completion in my general construction.

\subsection{The Torus}
Consider $M=\T^2$ symplectic, with symplectic area $2\pi/\hbar$. There are several ways to quantize this manifold. Traditional geometric quantization corresponds to taking the pair groupoid. This is prequantizable if and only if $1/\hbar$ is an integer, $m$. For any well-behaved polarization, this quantizes to the algebra of $\abs{m}\times\abs{m}$ matrices, which  isn't very interesting.

The more natural and interesting choice is the fundamental groupoid. There are several nice choices of polarization here; each produces the same algebra (the noncommutative torus) in significantly different ways. First, we could choose a \Kahler\ polarization on $\T^2$ and the corresponding exact polarization of $\PI(\T^2)$. This corresponds to the construction I studied in \cite{haw5} and described in the previous section.

We can always identify $\T^2$ as a quotient $\T^2\cong V/2\pi\,\Z^2$ of a symplectic vector space $V=\R^2$ with constant symplectic form $\omega = \tfrac{1}{2\pi\hbar} dx^1\wedge dx^2$. The fundamental groupoid is identified with $T^*\T^2\cong \T^2\times V^*$. We can employ a ``horizontal'' polarization $\P_{\mathrm{hor}}$, which is the kernel of the groupoid fibration $\T^2\times V^* \onto V^*$.  

The computation of the cocycle for a Poisson vector space (Sec.~\ref{Constant}) is still valid here, so we have $\sigma_0:V^*\times V^*\to\T$,
\[
\sigma_0(y,y') = e^{\pi i\hbar (y_2y'_1-y_1y'_2)} .
\]
The leaves of $\P_{\mathrm{hor}}$ are $2$-tori, so the holonomy around these gives two Bohr-Som\-mer\-feld conditions. 

Although the symplectic form on this groupoid is exact, it does not admit an adapted symplectic potential. This doesn't matter because the prequantization is unique and there is no danger of contradiction. The simplest choice of symplectic potential is 
\[
\theta = y_idx^i .
\]
For the fiber over $y\in V$, the holonomies for the two generators of $\pi_1(\T^2)$ are $e^{2\pi i y_1}$ and $e^{2\pi i y_2}$. The holonomy is therefore trivial if and only if $y\in\Z^2$. The reduced group(oid) is thus $\Z^2$ and the quantization is $\cs(\Z^2,\sigma_0)$. This \cs-algebra is generated by two unitaries $U$ and $V$ such that $UV = e^{-2\pi i\hbar}VU$. This is one of the standard presentations of the noncommutative torus algebra.

This construction generalizes easily to $\T^n$ with a constant Poisson structure. The quantization is $\cs(\Z^n,\sigma_0)$, where $\sigma_0$ is constructed by exponentiating the Poisson structure.

Another real polarization for $\T^2$ was considered by Weinstein \cite{wei3}. Let $\P_{\mathrm{Wei}}$ be the foliation whose leaves are the cylinders of constant $x^2$ and $y_1$. This is indeed a polarization, and the quotient groupoid is the action groupoid of $\R$ acting on $S^1$ by $(x^2,y_1)\mapsto x^2 - 2\pi \hbar y_1$.

This polarized symplectic groupoid is a quotient of $T^*\Pair(\R)$ with the vertical polarization. As such, the cocycle on $T^*\T^2/\P_{\mathrm{Wei}}$ is trivial, but there is one Bohr-Sommerfeld condition to take into account. The holonomy around the leaf over $(x^2,y_1)$ is $e^{2\pi i y_1}$. This is trivial if and only if $y_1\in\Z$. The reduced groupoid is thus the action groupoid of $\Z$ acting on $S^1$ by $2\pi \hbar$ rotations. The quantization is the crossed product \cs-algebra:
\[
\cs_{\P_{\mathrm{Wei}}}[\PI(\T^2),\sigma] \cong \cs(\Z\rtimes S^1) = \cs[\Z,\C(S^1)] .
\]
This is another standard presentation of a noncommutative torus algebra and gives it the name ``irrational rotation algebra'' (if $\hbar$ is irrational).

This general framework for quantization provides a heuristic explanation as to why these different constructions give the same algebras. It is simply a case in which the algebra is independent of the polarization.

\section{More About Polarizations}
\label{Polarization2}
\subsection{Real Polarizations of Lie Algebroids}
If $(\gamma,\eta)\in\G_2$ is a pair of composable elements in a Lie groupoid, and we happen to know the values of a polarization at $\gamma$ and $\eta$, then the definition of polarization tells us what the polarization is at $\gamma\eta$. So, if $U\subset \G$ is a subset that generates the entire groupoid, then any polarization can be reconstructed from its restriction to $U$. This suggests that for an $\s$-connected groupoid, any polarization can be reconstructed from its restriction to an infinitesimal neighborhood of the identity submanifold.

In other words, a polarization of $\G$ should be completely described by an infinitesimal structure in the same way that $\G$ is almost completely described by its Lie algebroid. It seems appropriate to call this hypothetical structure \emph{a polarization of a Lie algebroid}. 

Lie algebroids have many advantages over Lie groupoids. It is usually easier to explicitly construct a Lie algebroid than its groupoid. Simple Lie algebras are easier to classify than simple Lie groups. Polarization of Lie algebroids should have these same advantages. It may even be possible to carry out my quantization recipe in some cases without having to explicitly construct a symplectic groupoid.

The general description of Lie algebroid polarizations is a delicate matter, and I leave it to a future paper. However, the case of \emph{real} polarizations is much simpler.

For any Lie algebroid $A$ over $M$, the tangent space $TA$ is a \emph{double Lie algebroid} \cite{mac3}. That is, it has two structures as a Lie algebroid: over $A$ and over $TM$.
\begin{definition}
\label{A_polarization}
A \emph{sub double Lie algebroid} of a double Lie algebroid is a subset that is a Lie subalgebroid of both structures.
A \emph{real polarization of a Lie algebroid} $A$ is a sub double Lie algebroid $\p\subset TA$ that is wide over $A$.
\end{definition}

This definition is justified by the following result. Recall from Theorem~\ref{restatement} that a real polarization of a groupoid $\G$ is a sub Lie algebroid-groupoid $\P\subset T\G$. 
\begin{thm}
If $\P\subset T\G$ is a real polarization of a Lie groupoid, then $\Alg(\P)\subset T\Alg(\G)$ is a real polarization of the Lie algebroid $\Alg(\G)$.
\end{thm}
\begin{proof}
The functor $\Alg$ maps Lie algebroid-groupoids to double Lie algebroids and respects subalgebroids. See \cite{mac3} for the proof that $T\Alg(\G)=\Alg(T\G)$.

A polarization $\P\subset T\G$ is a wide Lie subalgebroid because it has a fiber over every point of $\G$. This implies that $\Alg(\P)$ is wide over $\Alg(\G)$.
\end{proof}

Definition \ref{A_polarization} is succinct, but it is useful to unpack it into more familiar geometrical structures.

To distinguish the two bundle structures of $TA$, I will refer to the tangent bundle structure as ``vertical'' and the other structure as ``horizontal''. The four different bundle structures involved can be summarized in the following (commutative) diagram:
\beq
\label{TA}
\begin{CD}
TA @>>> TM \\
@VVV @VVV \\
A @>>> M .\\
\end{CD}
\eeq

As a double Lie algebroid, $TA$ is in particular a double vector bundle.
One of the most important constructs from a double vector bundle is the \emph{core} \cite{mac3}. This is the intersection of the kernels of the vertical and horizontal projections. The core is an ordinary vector bundle over the base manifold $M$. In this case, the core of $TA$ is naturally identified with $A$.

The horizontal Lie algebroid structure of $TA$ is defined using two special types of sections. If $\xi\in\Gamma(M,A)$, then $T\xi\in\Gamma(TM,TA)$ is just the result of applying the tangent functor to $\xi:M\to A$. The core section $\hat\xi\in\Gamma(TM,TA)$ is defined by the identification of $A$ with the core of $TA$. These two types of sections span $\Gamma(TM,TA)$ as a $\C^\infty(TM)$-module, so the horizontal Lie algebroid structure is completely defined by:
\begin{gather*}
\#_{TA}T\xi = T(\#_A\xi) \mbox, \quad
\#_{TA}\hat\xi = \widehat{\#_A\xi}\mbox, \\
[T\xi,T\zeta] = T([\xi,\zeta]) \mbox, \quad
[T\xi,\hat\zeta] = \widehat{[\xi,\zeta]}\mbox, \quad \text{and}\quad
[\hat\xi,\hat\zeta] = 0
\end{gather*}
for $\xi,\zeta\in\Gamma(M,A)$.

\begin{definition}
If $\nabla$ is a flat partial connection on a quotient bundle $A/\P_1$, then a section of $A$ is \emph{$\nabla$-stable} if it is $\nabla$-constant modulo $\P_1$, i.e., it gives a $\nabla$-constant section of $A/\P_1$.
\end{definition}

\begin{thm}
\label{polarization_restatement}
A real polarization $\p\subset TA$ determines a triple $(\P_0,\P_1,\nabla)$ where $\P_0\subseteq TM$ is a foliation, $\P_1\subseteq A$ is a subbundle, and $\nabla$ is a flat $\P_0$-connection on $A/\P_1$, such that
over any open $U\subseteq M$:
\begin{enumerate}
\item
The $\nabla$-stable sections of $A$ form a Lie subalgebra of $\Gamma(U,A)$;
\item
the sections of $\P_1$ form a Lie ideal inside this Lie subalgebra; 
\item
$\#\P_1\subseteq \P_0\:$;
\item
the map $A/\P_1\to TM/\P_0$ induced by the anchor intertwines $\nabla$ with the $\P_0$-Bott connection (Def.~\ref{Bott}).
\end{enumerate}
Conversely, such a triple determines a real polarization of $A$, such that the structures are equivalent.
\end{thm}

\begin{proof}
Define $\P_0\subseteq TM$ to be the horizontal base of $\p$. Wide in the vertical structure means that the vertical base of $\p$ is $A$. This implies that $\P_0$  is a wide subalgebroid of $TM$ --- in other words, a foliation of $M$.
Graphically, this can be summarized by a subdiagram of \eqref{TA}:
\[
\begin{CD}
\p @>>> \P_0 \\
@VVV @VVV \\
A @>>> M .\\
\end{CD}
\]

The core of $TA$ is canonically isomorphic to $A$ itself. Define $\P_1\subseteq A$ to be the core of $\p\subset TA$.

By definition, $\p\subset TA$ is a sub double Lie algebroid. Compatibility with both vector bundle structures implies that this is actually a linear subalgebroid in both the vertical and horizontal structures, i.e.,  all relevant maps are fiber linear in the other bundle structure.

Being a linear subalgebroid in the vertical structure means that $\p$ is a linear foliation of $A$. (If $\P_0=TM$ and $\P_1=0$, then this would be a flat linear connection on $A$.) In general, a linear foliation of $A$ is a flat $\P_0$-connection on $A/\P_1$. Define $\nabla$ to be this connection.

Now consider the horizontal Lie algebroid structure. For $\xi\in\Gamma(U,A)$, $\hat\xi$ is a section of $\p$ if and only if $\xi$ is a section of $\P_1$; $T\xi$ is a section of $\p$ if and only if $\xi$ is $\nabla$-stable.

We can now translate the condition that $\p\subset TA$ be a horizontal subalgebroid into conditions on $\P_0$, $\P_1$, and $\nabla$. Let $\xi,\zeta\in\Gamma(U,A)$ be $\nabla$-stable sections and $\lambda\in\Gamma(U,\P_1)$, so that $T\xi,T\zeta,\hat\lambda \in \Gamma(TU,\p)$.

The brackets $[T\xi,T\zeta]=T([\xi,\zeta])$ and $[T\xi,\hat\lambda] = \widehat{[\xi,\lambda]}$ must be sections of $\p$. Therefore $[\xi,\zeta]$ is $\nabla$-stable and $[\xi,\lambda]$ is a section of $\P_0$. 

The anchor $\#_{TA}\hat\lambda$ must be a section of $T\P_0$, so $\#_A\lambda$ is a section of $\P_0$. The anchor of $T\xi$ must be a section of $T\P_0$, so $\#_A\xi$ is stable with respect to the Bott connection of the foliation $\P_0$.

Conversely, if $\nabla$ is a connection satisfying these conditions, then over a sufficiently small neighborhood of any point, there exist enough $\nabla$-stable sections to generate the module of horizontal sections of $\p$ in this way. This then shows that $\p\subset TA$ is a linear horizontal subalgebroid.
\end{proof}

Theorem \ref{polarization_restatement} is very close to the definition \cite{h-m1,mac3} of an ideal system for $A$. An ideal system is for a Lie algebroid what a smooth congruence (Def.~\ref{congruence}) is for a Lie groupoid; it determines a quotient algebroid. A surjective submersion $p:M\onto M'$ determines an equivalence relation by
\[
x\sim y \Longleftrightarrow p(x)=p(y)
\]
for $x,y\in M$. This equivalence relation is a Lie subgroupoid $M\times_p M \subset \Pair(M)$. 
An ideal system of $A$ consists of a surjective submersion $p:M\onto M'$, a subbundle $\P_1\subset A$, and an action of $M\times_p M$ on $A/\P_1$. The definition is otherwise the same as  the conclusion of Theorem~\ref{polarization_restatement}, with $M\times_p M$ playing the role of the foliation $\P_0$ and the action playing the role of $\nabla$.

An ideal system always determines a real Lie algebroid polarization. If the fibers of $p$ are connected, then the ideal system is completely equivalent to this polarization. A real Lie algebroid polarization determines an ideal system if and only if $\P_0$ is simple (the kernel foliation of a submersion) and $\nabla$ has trivial holonomy.

To recover a groupoid polarization from a Lie algebroid polarization, one should apply the groupoid integration functor $\Grp$ to the horizontal structure of $\p$. In principle, there should be an integration theorem for real polarizations, saying that under some global conditions an algebroid polarization integrates to a groupoid polarization. However, I will not dwell on this here. It is a special case of the problem of integrating complex polarizations, and I plan to explore that in a future paper.

When $\p=\Alg(\P)$, the structures $\P_0$, $\P_1$ and $\nabla$ can be understood geometrically. First, $\P_0$ is the horizontal base of $\Alg(\P)$, so it is the unit manifold of $\P$ (as a groupoid); in other words, $\P_0 = \P \cap TM$. Thought of as foliations, $\P_0$ is the restriction of $\P$ to $M\subset \G$. 

As a vector bundle, $A=\Alg(\G)$ is defined to be the normal bundle to $M\subset \G$. This can be stated as an exact sequence of vector bundles, 
\beq
\label{A.sequence}
0\to TM \xrightarrow{T\unit} \unit^*T\G \to A \to 0.
\eeq
The subbundle $\P_1$ is given by the exact subsequence, 
\beq
\label{P.sequence}
0 \to \P_0 \to \unit^*\P \to \P_1 \to 0 .
\eeq
There are three natural ways of splitting \eqref{A.sequence}: using $\s$, $\t$, or $\inv$. The last is the most symmetrical; it identifies $\unit^*T\G$ with $TM\oplus A$ as the $\pm1$ eigenspace decomposition of $T\inv|_M$.
Because $\P$ is a real polarization, $T\inv|_M$ is an automorphism of $\unit^*\P$, so this is also a splitting of \eqref{P.sequence}.

The $\P$-Bott connection on $\P^\perp$  restricts along $M$ to a flat $\P_0$-connection on $\unit^*\P^\perp$. The foliation $\P$ is preserved by $\inv:\G\to\G$, so this connection must be compatible with $\inv$ and respect the splitting. The connection decomposes into the $\P_0$-Bott connection on $\P_0^\perp$ and the connection $\nabla$ on $\P_1^\perp$ (or equivalently, on its dual $A/\P_1$).

The case of a complex polarization is much more difficult to study, because $\P\subset T_\co\G$ is no longer a subgroupoid, it may be transverse to $M\subset \G$, and $T\inv|_M$ is no longer an automorphism of $\unit^*\P$.

Over any $x\in M$, $\ker \#_x\subseteq A_x$ is a Lie algebra (quite possibly $0$). If $\xi,\zeta\in\Gamma(M,A)$ satisfy $\#\xi(x)=\#\zeta(x)=0$, then the Lie algebra bracket $[\xi(x),\zeta(x)]$ is (by definition) the value of the Lie algebroid bracket $[\xi,\zeta]$ at $x$.
\begin{thm}
\label{P.ideal}
For any real polarization of $A$ and any point $x\in M$,
\[
\P_1\cap \ker \#_x \subseteq \ker \#_x
\]
is a Lie algebra ideal.
\end{thm}
\begin{proof}
Let $U\ni x$ be a neighborhood over which $\P_0$ is a simple foliation. Then any pair of vectors in $\P_1\cap\ker\#_x$ and $\ker\#_x$ are the values at $x$ of a section $\xi\in\Gamma(U,\P_1)$ and a $\nabla$-stable section $\zeta\in\Gamma(U,A)$. By the properties of a real polarization, $[\xi,\zeta] \in \Gamma(U,\P_1)$ and
\[
[\xi(x),\zeta(x)] = [\xi,\zeta](x) \in \P_1\cap\ker\#_x .
\]
\end{proof}

\begin{example}
One extreme case would be a polarization with $\P_0=TM$ and $\P_1=0$. Then $\nabla$ is a flat connection on $A$. If $\nabla$ has no holonomy, then it gives a trivialization of $A$ as a bundle. The constant sections form a Lie subalgebra $\mathfrak g \subset \Gamma(M,A)$. The anchor map gives an action of $\mathfrak g$ on $M$. In this case, $A$ must be the action Lie algebroid of $\mathfrak g$ on $M$.
\end{example}

\subsection{Real Polarizations of Poisson Manifolds}
Now consider the case that $\P$ is a real polarization of a symplectic groupoid $\Sigma$ integrating a Poisson manifold $M$. In this case, $\Alg(\Sigma) = T^*M$, so what are the properties of $\Alg(\P) \subset T\espace T^*M$?

We can (again) use $\inv$ to split the exact sequence \eqref{A.sequence} and identify $\unit^*T\Sigma$ with $TM\oplus T^*M$. In this identification, the symplectic structure on $\unit^*T\Sigma$ is given by the pairing between $TM$ and $T^*M$. Because $\P$ is a symplectic polarization, $\unit^*\P \cong \P_0\oplus \P_1 \subset TM\oplus T^*M$ must be Lagrangian. This just means that $\P_1=\P_0^\perp$. The $\P_0$-connection $\nabla$ thus acts on $\P_0^{\perp\perp}=\P_0$ itself.

For any isotropic foliation of a symplectic manifold, the Bott connection gives (through the symplectic form) a partial connection on the foliation itself --- and this is torsion-free (see \cite{vai2}). The connection coming from $\Alg(\P)$ is just the restriction of this connection, so it is also torsion-free.

\begin{definition}
\label{PP.def}
A \emph{real polarization of a Poisson manifold} is a pair $(\P_0,\nabla)$ where $\P_0\subseteq TM$ is a coisotropic foliation of $M$ and $\nabla$ is a flat, torsion-free $\P_0$-connection such that over any open $U\subset M$:
\begin{enumerate}
\item
The set of $\nabla$-stable $1$-forms is closed under the Koszul bracket;
\item
the sections of $\P_0^\perp$ form a Lie ideal inside this;
\item
if $f\in\C^\infty(U)$ is $\P_0$-constant then $\nabla \#df =0$.
\end{enumerate}
\end{definition}
``Coisotropic'' is just a restatement of $\#\P_1 = \#\P_0^\perp \subseteq\P_0$. Everything in this definition is a restatement of the properties of a real Lie algebroid polarization to this case, except the torsion condition. The vanishing torsion of $\nabla$ is necessary to insure that $\P$ remains Lagrangian away from $M\subset \Sigma$, as illustrated by this example.
\begin{example}
Consider $S^3$ with $0$ Poisson structure.  This is integrated by $\Sigma=T^*S^3$ with the structure of a bundle of Abelian groups. Identifying $S^3\cong \mathrm{SU}(2)$, the right-invariant trivialization of $TS^3$  defines a linear map $p:T^*S^3\to \R^3$ which is a groupoid homomorphism and in fact a fibration. So, its kernel foliation $\ker Tp$ is a groupoid polarization. 

The Lie algebroid polarization $\Alg(\P)$ is given by $\P_0=TS^3$, $\P_1=0$, and the right invariant connection. 
This connection is not torsion-free, so this does not satisfy Definition~\ref{PP.def}.

The map $p$ would be Poisson if we gave $\R^3$ the nonzero Lie-Poisson structure of $\mathfrak{su}(2)^*$, therefore $\ker Tp$ is not Lagrangian and is not a \emph{symplectic} groupoid polarization.
\end{example}

What if $M$ is symplectic? A symplectic polarization does not require the additional structure of a connection, so these two definitions of polarization need to be reconciled. 
\begin{example}
If $M$ is symplectic and $\P_0$ is a Lagrangian foliation, then the last condition in Definition \ref{PP.def} completely fixes a connection $\nabla$. Any vector in $\P_0$ is a value of a Hamiltonian vector field $\#df$ for a $\P_0$-constant function $f$ defined over some open set. Any such vector field is $\nabla$-constant. So, the connection isn't really an additional structure in this case, and a real symplectic polarization is a special case of a real Poisson polarization.
\end{example}
However, depending upon the topology of a symplectic manifold, there may exist real Poisson polarizations for which $\P_0$ is not Lagrangian and $\nabla$ does carry additional information.

\begin{example}
In Section \ref{Constant}, I described the quantization of vector space $V$ with constant Poisson structure, using a ``horizontal'' polarization. In this case, $\P_0=TV$ (which is automatically coisotropic) and $\nabla$ is the trivial connection of the vector space.
\end{example}

At a point $x\in M$ of a Poisson manifold, the Lie algebra $\ker \#_x$ has a particular significance. It is the dual of the linearization of the Poisson structure transverse to the symplectic leaf through $x$; see \cite{wei0}. This is only nonabelian if $x$ is a singular point of the Poisson structure. 

Theorem \ref{P.ideal} shows that $\P_0^\perp\cap\ker\#_x$ is a Lie algebra ideal in $\ker\#_x$. This can be a very restrictive condition. 
\begin{prop}
Let $M$ be a connected Poisson manifold admitting a real polarization whose connection is geodesically complete. If the Poisson structure vanishes at $x\in M$ and the linearization is the dual of a \emph{simple} Lie algebra $\mathfrak g$, then $M$ is isomorphic to $\mathfrak g^*$ with the Lie-Poisson structure.
\end{prop}
\begin{proof}
Since $\mathfrak g$ has no nontrivial ideals, $\P_{0\,x}$ must be all or nothing. By local triviality, we must have either $\P_0=0$ or $\P_0=TM$, but $0$  is not coisotropic. So, $\P_0=TM$ and $\nabla$ gives $M$ a locally affine structure. 

Since $\nabla$ is geodesically complete, the universal cover $\Mt$ is an affine space, but we can assign a preimage of $x$ to be $0$, thus making $\Mt$ a vector space. Because $\P_0^\perp=0$, $\nabla$-stable just means $\nabla$-constant. A function $f\in\C^\infty(\Mt)$ is affine if and only if $0=\nabla df$. So, for two affine functions $f,g\in\C^\infty(\Mt)$,
\[
0 = \nabla [df,dg]_\pi = \nabla d\{f,g\} .
\]
Hence $\{f,g\}$ is affine. Since the Poisson structure vanishes at $0$, the bracket $\{f,g\}$ vanishes at $0$ and is thus linear. So, the linear functions form a Lie algebra which must coincide with the linearization $\mathfrak g$.
This identifies $\Mt$ as $\mathfrak g^*$. Because $\mathfrak g$ is simple, the Lie-Poisson structure only vanishes at one point, and thus we must have $M\cong \Mt\cong \mathfrak g^*$
\end{proof}

The group $\mathrm{SU}(2)$ admits an essentially unique Poisson-Lie group structure, corresponding to the standard quantum groups $\mathrm{SU}_q(2)$. (This should not be confused with the ``Lie-Poisson'' structure on the dual of a Lie algebra.)
\begin{prop}
There does not exist any real polarization of this Poisson structure on $S^3\cong \mathrm{SU}(2)$.
\end{prop}
\begin{proof}
The Poisson bivector vanishes along a single great circle $C\subset S^3$. The nontrivial symplectic leaves are open $2$-dimensional discs with boundary $C$.

$S^3$ does not admit a flat, torsion-free connection, so a real polarization must have $\rk \P_0=1$ or $2$.

At any point of $C$, the linearized Lie algebra is isomorphic to that generated by $\{X,Y,Z\}$ with the brackets:
\[
[X,Y]=0,\quad [Y,Z]=X,\quad [X,Z]=-Y .
\]
The only proper ideal is spanned by $X$ and $Y$. Its annihilator is $1$-dimensional and corresponds to the tangent to $C$.

So, if there is a real polarization, then $\rk \P_0=1$, and $\P_0$ is tangent to $C$. Being coisotropic implies that $\P_0$ is tangent to the symplectic leaves. 

Restricting $\P_0$ to the closure of a symplectic leaf, we have a closed $2$-disc with a rank $1$ foliation that is tangential to the boundary. This is impossible.
\end{proof}

This type of problem can sometimes be circumvented by using complex polarizations. Theorem~\ref{P.ideal} generalizes to complex polarizations, but the generalization of $\P_1$ is not locally trivial, and a Lie algebra sometimes has more ideals when it is complexified.

The structure of a real Poisson polarization is quite restrictive.
\begin{thm}
\label{real.structure}
Around any point of a real polarized Poisson manifold there exist coordinates $(x^i,y^\alpha)$ such that:
\begin{enumerate}
\item
The $\P_0$-leaves are the subspaces of fixed $y$;
\item
$\nabla_i = \frac{\partial}{\partial x^i}$;
\item
the $y$-$y$-components of $\pi$ vanish;
\item
the $x$-$y$-components of $\pi$ are constant in the $x$'s;
\item
the $x$-$x$-components of $\pi$ are affine in the $x$'s.
\end{enumerate}
\end{thm}
\begin{proof}
The flat, torsion-free $\P_0$-connection $\nabla$ gives a locally affine structure to the leaves of $\P_0$. 

Over some open $U\subset M$, let $f,g\in\C^\infty(U)$ be functions which are affine along the $\P_0$-leaves. Equivalently, $df$ and $dg$ are $\nabla$-stable. Definition \ref{PP.def} then requires that the Koszul bracket,
\[
[df,dg]_\pi = d\{f,g\}
\mbox,
\]
is also $\nabla$-stable. Therefore $\{f,g\}$ is also affine along the $\P_0$-leaves.

Let $h\in\C^\infty(U)$ be constant along the $\P_0$-leaves. This means that $dh$ is a section of $\P_0^\perp$. Definition \ref{PP.def} requires that the Koszul bracket,
\[
[df,dh]_\pi = d\{f,h\}
\mbox,
\]
is also a section of $\P_0^\perp$. Therefore, $\{f,h\}$ is also $\P_0$-constant.

Because of the locally affine structure, we can construct a foliated coordinate chart around any point of $M$ such that $\nabla$ is just given by partial derivatives. Let $x^i$ be the leafwise coordinates and $y^\alpha$ the transverse coordinates in such a chart.

The $y$'s are $\P_0$-constant, so $\pi^{\alpha\beta} = \{y^\alpha,y^\beta\} = 0$, because $\P_0$ is coisotropic.

The $x$'s are affine in the $\P_0$-leaves, so $\pi^{i\alpha} = \{x^i,y^\alpha\}$ is $\P_0$-constant.

Finally, $\pi^{ij} = \{x^i,x^j\}$ is affine in the $\P_0$-leaves, so it is affine in the $x$'s.
\end{proof}

\subsection{Computing the Twist}
The description of real polarizations at the Poisson manifold level opens the possibility of constructing the algebra $\cs_\P[\Sig(M),\sigma]$ without having to work with the symplectic groupoid $\Sig(M)$ directly.

Unfortunately, when Bohr-Sommerfeld conditions come into play, $\Sigma_\BS$ is typically not $\s$-connected. In that case, the Lie algebroid does not contain essential information about the groupoid. However, if the Bohr-Sommerfeld conditions are trivial, then there is still hope.

If $\P$ is a strongly admissible real polarization of a groupoid $\G$, then $\Alg(\P)$ is equivalent to an ideal system for $\Alg(\G)$. We can compute $\Alg(\G/\P)$ directly. Its base manifold is the leaf space $M/\P_0$. The Lie algebra of sections is identified with the space of $\nabla$-stable sections of $\Alg(\G)$ modulo the sections of $\P_0$.

The description of the Lie algebroid cohomology complex of $\Alg(\G/\P)$ is actually simpler. It consists of the $\nabla$-constant sections of $\Wedge^\bullet\P_1^\perp$.

Now, let $M$ be a Poisson manifold and $\P$ a strongly admissible real polarization of $\Sig(M)$ with simply connected leaves. Let $p:\Sig(M)\to\Sig(M)/\P$ be the quotient map. Let $A:= \Alg[\Sig(M)/\P]$ be the quotient Lie algebroid.

Prequantization of $\Sig(M)$ gives a cocycle $\sigma$ which has a class $[\sigma]\in\Tw[\Sig(M)]$. This should be equivalent to the pull-back of the class of the reduced cocycle $[\sigma_0]\in\Tw(\Sig(M)/\P)$. Applying the ``characteristic class'' map $\Psi$ from Section~\ref{Prequantization} gives a commutative diagram,
\[
\begin{CD}
\Tw[\Sig(M)] @>{\Psi}>> H^2_\pi(M) \\
@A{\Tw(p)}AA @A{p^*}AA \\
\Tw[\Sig(M)/\P] @>{\Psi}>> H^2_{\mathrm{Lie}}(A) .
\end{CD}
\]
We know that $\Psi(\sigma) = [\pi]\in H^2_\pi(M)$, so the pull-back of $\Psi(\sigma_0)$ must be cohomologous to $\pi$. 

If $p^*:H^2_{\mathrm{Lie}}(A) \to H^2_\pi(M)$ happens to be injective, then in light of Theorem~\ref{Psi_injective}, the problem comes down to finding a $\nabla$-constant bivector $c\in\Gamma(M,\Wedge^2\P_0)$ such that
\[
c = \pi + \delta X
\]
where $\delta X$ is the Poisson differential of some vector field. We can then recover $[\sigma_0]$ from $[c]=\Psi(\sigma_0)$, and only the class of $\sigma_0$ is relevant to computing the algebra.

It is not obvious the question of whether such a $c$ exists, but it does in the very general example that I will discuss in Section~\ref{Affine}.

\subsection{Totally Complex Polarizations}
In the quantization of symplectic manifolds, the best kind of polarization is a totally complex polarization given by a positive \Kahler\ structure. In that case, all the elegant tools of complex analysis and Dirac-type operators apply. One might expect totally complex polarizations to be equally useful for symplectic groupoids, but one should be prepared for disappointment, because such polarizations don't exist for most groupoids.

\begin{thm}
\label{totally_complex}
If $\G$ is an $\s$-connected Lie groupoid over $M$ with a complex structure such that the antiholomorphic tangent bundle is a groupoid polarization, then $\G$ is a covering of $\Pair(M)$, the anchor is an isomorphism $\#:\Alg(\G)\isom TM$ of Lie algebroids, and there exists a unique complex structure on $M$ that induces the complex structure on $\G$. 
\end{thm}
\begin{proof}
Again, $\unit^*T\G$ can be identified with $TM\oplus \Alg(\G)$ using the $\pm1$ eigenspace decomposition by $T\inv$. Each fiber of this is a vector space groupoid, and so the product can be written in terms of the anchor (as in Sec.~\ref{Constant}).

Let $J:T\G\to T\G$ be the complex structure, and $\P\subset T_\co\G$ the antiholomorphic tangent bundle.
By Hermiticity, $T\inv:\P\to\bar\P$ (and \emph{vice versa}), so $J\circ T\inv = -T\inv\circ J$. The complex structure intertwines the eigenvalues of $T\inv$, so it restricts to a vector bundle isomorphism $J:TM\isom\Alg(\G)$. The antiholomorphic tangent bundle along $M$ is
\[
\unit^*\P = \{(v,iJv)\mid v\in T_\co M\} .
\]
So, take an arbitrary $v\in TM$. 
By Multiplicativity, $(v,iJv)$ can be written as a product,
\begin{align*}
(v,iJv)&=(u,iJu)\cdot (w,iJw)\\
 &= (u-\tfrac{i}2\#Jw,iJ[u+w]) = (w+\tfrac{i}2\#Ju,iJ[u+w]) .
\end{align*}
The two expressions occur because $(u,iJu)$ and $(w,iJw)$ must be composable.
The first implication of this is that $v=u+w$. Applying $T\s$ to this gives the  equations,
\[
T\s(w,iJw)=  T\s(v,iJv) = T\s(u+w,iJ[u+w]).
\]
The map $T\s:\unit^*T\G\to TM$ is linear, so
\[
0 = T\s(u,iJu) = u - \tfrac{i}2\#J u .
\]
Likewise, applying $T\t$ gives,
\[
0 = T\t(w,iJw) = w + \tfrac{i}2\#J w .
\]
This shows that
\[
\tfrac12 \#Jv = i(-u+w)
\]
and
\[
\tfrac12\#J \tfrac12\#J v = -v \mbox,
\]
which means that $J_0:= \tfrac12\#J$ is an almost complex structure and $\#$ must be an isomorphism.  

Since $\#:\Alg(\G)\to TM$ is always a Lie algebroid homomorphism, it is now a Lie algebroid isomorphism. This makes the groupoid anchor $(\t,\s):\G\to \Pair(M)$ a local homoeomorphism.

From the explicit presentation of $\unit^*\P$, we can compute
\[
\P_0 = T\t(\unit^*\P) = \{v+iJ_0v\mid v\in T_\co M\} = \{v\in T_\co M\mid J_0v=-iv\} .
\]
By Lemma \ref{P0.lem}, $\P_0$ is involutive, so $J_0$ is integrable, and $\P_0$ is the antiholomorphic tangent bundle for this complex structure.

By Lemma \ref{P0.lem}, for any $\gamma\in\G$, $T\t(\P_\gamma) = \P_{0\,\t(\gamma)}$ and
\[
T\s(\P_\gamma) = T\t\circ T\inv(\P_\gamma) = T\t(\bar\P_{\gamma^{-1}}) = \bar \P_{0\,\s(\gamma)}
\mbox,
\]
so
\[
\P \subseteq T\t^{-1}\P_0 \cap T\s^{-1}\bar\P_0 = T(\t,\s)^{-1}(\P_0\times\bar\P_0) .
\]
This is an equality along the unit manifold. Computing the rank shows that it is an equality everywhere. 
\end{proof}

\begin{cor}
If $\Sigma$ is an $\s$-connected symplectic groupoid with a totally complex polarization, then its base manifold is a \Kahler\ manifold and the polarization is the corresponding exact polarization.
\end{cor}
 
Despite this negative result, there is a nice class of polarizations associated to complex structures on Poisson manifolds; I describe this in Section~\ref{Kaehler-Poisson}.

Note that Hermiticity was essential to Theorem \ref{totally_complex}.
\begin{example}
If $G$ is a complex algebraic group, then the (anti)holomorphic tangent bundle is \emph{not} a polarization. It is involutive and multiplicative, but not Hermitian.
\end{example}

\subsection{Induced Polarizations}
\label{induced}
Here, I give two techniques by which a polarization on one groupoid can determine a polarization on another. 

First, we need a few definitions. A Lie algebra-groupoid is just a special case of a Lie algebroid-groupoid (see Def.~\ref{LA-groupoid}); this is also known as a strict 2-Lie algebra.
\begin{definition}
If $\h$ is a Lie algebra-groupoid, $\G$ is a Lie groupoid, and $\daction :\h\to\X^1(\G)$ is a Lie algebra action, then $\daction$ is \emph{multiplicative} if $\daction:\h\times\G\to T\G$ is a groupoid  homomorphism. A subspace $F\subset \h_\co$ of the complexification is \emph{multiplicative} if $F\cdot F = F$ and \emph{Hermitian} if $\inv(F)=\bar F$.
\end{definition}
\begin{lem}
\label{push_forward}
Let $\daction:\h\to \G$ be a multiplicative action of a Lie algebra-groupoid on a groupoid, and $F\subset \h_\co$ a multiplicative, Hermitian Lie subalgebra. If
\[
\Im \daction(F) := \{\daction_v(\gamma) \mid v\in F, \gamma\in\G\} \subset T_\co\G
\]
is a regular distribution, then it is a polarization of $\G$.
\end{lem}
\begin{proof}
In particular, $\daction : F \to \X^1(\G)$ is a Lie algebra action, therefore $\Im \daction(F)$ is involutive.

Multiplicativity of $\daction$ and Hermiticity of $F$ imply that
\[
\inv_{T\G}\circ \daction(F) = \daction\circ \inv_\h(F) = \daction(\bar F) .
\]
So, $\Im \daction(F)$ is Hermitian.

Finally, multiplicativity of $\daction$ and $F$ implies multiplicativity of $\Im \daction(F)$.
\end{proof}
In particular, there is the canonical action (by right-invariant vector fields) of $\h$ on the Lie group-groupoid (strict 2-group) $\Grp(\h)$. $F$ determines a right-invariant polarization of $\Grp(\h)$. Any right-invariant polarization is of this form, so a multiplicative, Hermitian Lie subalgebra $F\subset \h_\co$ is equivalent to a right-invariant polarization of $\Grp(\h)$.

\begin{lem}
\label{inverse.image}
Let $q:\G\to\G'$ be a groupoid homomorphism, and $\Q$ a polarization of $\G'$. If the inverse image
\[
\P:= Tq^{-1}(\Q) \subset T_\co\G
\]
is a regular distribution, then it is a polarization of $\G$. 
\end{lem}
\begin{proof}
If $X,Y\in\Gamma(\G,\P)$ are $q$-projectable, then because $\Q$ is involutive, 
\[
Tq([X,Y]) = [Tq(X),Tq(Y)]
\]
 is a section of $\Q$. Therefore $[X,Y]$ is a section of $\P$, and $\P$ is involutive. 

Hermiticity of $\Q$ implies Hermiticity of $\P$, because $q$ intertwines the inverses. 

Now use the multiplicativity of $\Q$.
If $(X,Y)\in \P_2$ is a pair of composable vectors from $\P$, then $Tq(X\cdot Y) = Tq(X)\cdot Tq(Y) \in \Q$, therefore $X\cdot Y\in \P$. 

If $(\gamma,\eta)\in\G_2$ are composable and $Z\in \P_{\gamma\eta}$, then there exist $X'\in \Q_{q(\gamma)}$ and $Y'\in \Q_{q(\eta)}$ such that $X'\cdot Y' = Tq(Z)$. Choose an arbitrary vector $X \in \P_\gamma$ such that $Tq(X) = X'$. If we \emph{define} $Y$ by $X\cdot Y = Z$, then
\[
Tq(Y) = Tq(X^{-1}\cdot Z) = X'^{-1}\cdot Tq(Z) = Y' \in \Q
\]
so $Y\in \P_\eta$. This proves multiplicativity of $\P$.
\end{proof}
The simplest example is when $q$ is a fibration and $\Q=0$; then $Tq^{-1}(0) = \ker Tq$ is the kernel foliation. A little more generally, if $q$ is a fibration and $\Q$ is a strongly admissible real polarization, then $Tq^{-1}(\Q)$ is the kernel foliation of the composition of $q$ with the quotient fibration $\G'\onto\G/\Q$.

\section{Further Examples}
\label{Examples2}
\subsection{Trivial Poisson Structure}
If the Poisson structure vanishes, then quantization should return the algebra $\C_0(M)$. This is very important for a consistent classical limit.

First, suppose that we start with a manifold $M$ with Poisson bivector $\pi$ and symplectic integration $\Sigma$. Choose some polarization $\P$ of $\Sigma$. If the Poisson bivector is rescaled to $\hbar\pi$, then the symplectic form of $\Sigma$ is rescaled to $\omega/\hbar$.

Now, imagine taking the classical limit $\hbar\to 0$. In doing this, we should ``zoom in'' around the unit manifold of $\Sigma$ to balance the growing symplectic form. If we do so, then in the limit we are left with $T^*M$ with the canonical symplectic form and addition on fibers as the groupoid operation. This is the unique $\s$-simply connected integration of $M$ with $0$ Poisson structure.

In this process, the polarization becomes stretched vertically. What is left is the direct sum of the horizontal part $\P\cap T_\co M$ and its annihilator. This is (the complexification of) a real polarization of $T^*M$, although it may have developed some singularities.

So, \emph{all polarizations become real in the classical limit}.  
A (regular) real Poisson polarization of $M$ with $\pi=0$ is just a foliation $\P_0$ of $M$ with a locally affine structure on the leaves. 

With this in mind, consider a manifold $M$ with trivial Poisson bivector $\pi=0$, symplectic groupoid $T^*M$ and a polarization $\P$.

The classical limit of a \Kahler\ polarization is a vertical polarization, i.e., $\P_0=0$.
\begin{example}
If $\P_0=0$, then the fibers of $\P$ are the cotangent fibers. These are simply connected, so there are no Bohr-Sommerfeld conditions. The Liouville form is a multiplicative adapted symplectic potential, so the reduced cocycle $\sigma_0$ is trivial. The reduced groupoid is the quotient $T^*M/\P = M$ with trivial groupoid structure. Therefore the quantization recipe gives
\[
\cs_\P(T^*M,\sigma) = \cs(M)= \C_0(M) .
\]
\end{example}

Any other (not totally complex) polarization will have some horizontal part in the classical limit. Two elementary examples illustrate how this may still give the same algebra.
\begin{example}
Let $M=\R$. A Poisson polarization is given by $\P_0=T\R$ and the trivial connection $\nabla$. The symplectic groupoid polarization is the horizontal foliation of $T^*\R$. 

These leaves are simply connected, so there is no Bohr-Sommerfeld condition. The reduced groupoid is the quotient, $T^*\R/\P \cong \R$ as an additive group. The group cohomology of $\R$ is trivial, so the reduced cocycle $\sigma_0$ must be trivial. Therefore  
\[
\cs_\P(T^*\R,\sigma) = \cs(\R) \cong \C_0(\R) .
\]
\end{example}

This becomes more subtle when a Bohr-Sommerfeld condition comes into play, but it also illustrates the importance of the Bohr-Sommerfeld conditions.
\begin{example}
Let $M = S^1 = \R/2\pi\,\Z$. Again let $\P_0=TS^1$ with the trivial connection. Parametrize $T^*S^1$ with horizontal coordinate $x$ (modulo $2\pi$) and vertical coordinate $y$. The leaves of $\P$ are the circles of constant $y$, so the quotient groupoid is again the group $\R$.

The symplectic form is $dx\wedge dy$. 
Prequantization of $T^*S^1$ as a symplectic manifold is not unique, but since it is a symplectic groupoid, the holonomy around the units needs to be trivial. This makes the prequantization unique up to isomorphism as dictated by Theorem~\ref{c_prequantizability}. 

The holonomy around the $\P$-leaf at $y$ is $e^{-2\pi i y}$. So, the Bohr-Sommerfeld condition is $y\in \Z$. The reduced groupoid is the set of leaves with trivial holonomy; this is the additive group $\Z$.
Therefore
\[
\cs_\P(T^*S^1,\sigma) = \cs(\Z) \cong \C(S^1) .
\]
\end{example}

\label{trivial}
As a more general case, suppose that $\Action:\T^k\times M\to M$ is a free action of a torus $\T^k$ by diffeomorphisms of $M$. The orbits define a foliation $\P_0$ and a flat structure on the leaves; this is thus a Poisson polarization.

The quotient $M/\P_0=M/\T^k$ is a manifold; call it $N$. The original manifold $M$ is a principal $\T^k$-bundle over $N$.

The symplectic groupoid is again $T^*M$. The vertical part of the polarization is $\P_1=\P_0^\perp$, so $T^*M/\P_1\cong\P_0^*\cong \R^k\times M$. The quotient groupoid is thus the trivial bundle of groups $\R^k\times N$. The projection $q: T^*M\to\R^k\times N$ is the product of the dual map $T\espace\Action^*:T^*M\to\R^k$ and the composed projection $T^*M\to M\to N$.

The leaves of this polarization each have topology $\R^k\times\T^k$, so there are $k$ Bohr-Sommerfeld conditions. 
Suppose that $v$ is one of the integral basis vectors of the Lie algebra of $\T^k$. The holonomy around the flow of $v$ should be trivial. This holonomy can be computed using the Liouville form $\theta\in\Omega^1(T^*M)$. At $\xi\in T^*M$, the inner product of $T\espace\Action(v)$ with $\theta(\xi)$ is the pairing of $v$ with the dual map,
\[
T\espace\Action(v)\inner\theta(\xi) = \langle v, T\espace\Action^*(\xi)\rangle .
\]
The holonomy around the $v$ orbit through $\xi$ is given by exponentiating,
\[
1 \stackrel?= e^{2\pi i  \langle v, T\espace\Action^*(\xi)\rangle} .
\]
So the Bohr-Sommerfeld conditions require $T\espace\Action^*(\xi)\in \Z^k$. The reduced groupoid is
\[
T^*M_\BS/\P \cong \Z^k\times N .
\]

Let $w\in\Z^k$ be an integral vector. What is the reduced line bundle over $\{w\}\times N$? The space of sections over $\{w\}\times N$ is equal to the space of $\nabla = d+ i\theta$ constant functions over $q^{-1}(\{w\}\times N) = (T\espace\Action^*)^{-1}(w)$. In particular, these functions are constant along what is left of the cotangent fibers. If $\psi$ is such a function, and $v\in\R^k$ is a vector in the Lie algebra of $\T^k$, then 
\[
\Lie_{T\espace\Action(v)}\psi = -i\langle T\espace\Action(v),\theta\rangle \psi = -i\langle v,w\rangle \psi ,
\]
and the reduced line bundle over $\{w\}\times M$ is
\[
\co\times_{e^{-iw}} M 
\]
where $e^{-iw}$ is the equivalent unitary character (representation) of $\T^k$.

The twisted polarized convolution algebra is the tensor algebra over $\C^\infty_{\mathrm c}(N)$ generated by compactly supported smooth sections of these line bundles. This is isomorphic to a dense subalgebra of smooth functions on $M$. Effectively, what has happened is the decomposition of functions on $M$ into their Fourier components with respect to the $\T^k$ action. If this torus bundle is nontrivial, then the Fourier components are sections of these nontrivial line bundles.

\begin{example}
If $k=1$ then $M$ is a $\T$-bundle over $N$. The algebra of continuous functions on $M$ is generated (as a \cs-algebra) by sections of the associated line bundle over $N$.
\end{example}

\subsection{Bundle Affine}
\label{Affine}
A real polarization of a Poisson manifold is in particular a foliation with a locally affine structure on the leaves. Suppose that the leaves are actually the fibers of a vector bundle. Let $N:=M/\P_0$ be the base manifold and $q:M\to N$ the projection.

By Theorem~\ref{real.structure}, the space $\C^\infty_{\mathrm{aff}}(M)$ of fiber-affine functions is closed as a Lie algebra under the Poisson bracket. The space $q^*\C^\infty(N)\subset \C^\infty_{\mathrm{aff}}(M)$ is an Abelian Lie ideal. The quotient Lie algebra is naturally identified with the space of (homogeneous) fiber-linear functions or the space of sections of the dual vector bundle, which I will suggestively denote as $A\to N$. This gives a central extension of Lie algebras,
\[
0 \to \C^\infty(N) \stackrel{q^*}{\to} \C^\infty_{\mathrm{aff}}(M) \to \Gamma(N,A) \to 0 .
\]

Define $c\in\Gamma(N,\Wedge^2 A^*)$ to be the vertical part of the Poisson bivector along the $0$ section of $M=A^*$. Extending this to a fiber-constant bivector, we can decompose $\pi = c + \pi_0$.

The Lie bracket on $\Gamma(N,A)$ is given by $\pi_0$, so it is bidifferential, therefore $A$ is a Lie algebroid and $\pi_0$ is a Poisson structure. The Jacobi identity for $\pi$ implies that $c\in\Gamma(N,\Wedge^2A^*)$ is a Lie algebroid cocycle.

Now suppose that $(M,\pi)$ is integrable and we follow the quantization recipe using the groupoid $\Sig(M)$. The algebroid $A$ is the quotient of $T^*M$ by the polarization; the quotient of the groupoid is a fibration $p:\Sig(M)\to \G$ with connected fibers, therefore $\G=\Grp(A)$. To compute the quantization of $M$, we need the reduced cocycle $\sigma_0$ on $\G$, or rather its class $[\sigma_0]\in\Tw(\G)$.

By Theorem \ref{Psi_injective}, $[\sigma_0]$ is determined by its characteristic class $\Psi(\sigma_0)\in H^2_{\mathrm{Lie}}(A)$. 
This in turn is (partly) determined by $p^*(\Psi[\sigma_0]) = [\pi]\in H^2_\pi(M)$.

The pull-back map $p^*:H^\bullet(A)\to H^\bullet_\pi(M)$ is simple to describe. An algebroid cochain in $\Gamma(N,\Wedge^\bullet A^*)$ maps to the equivalent vertical, fiber-constant multivector on $M$. So, the class $p^*(\Psi[\sigma_0])$ includes a vertical, fiber-constant bivector on $M$ which is cohomologous to $\pi$.

Let $X\in\X^1(M)$ be the Euler vector field on $M$ (as a vector bundle over $N$). The eigenvectors of the Lie derivative operator $\Lie_X$ are the tensors which are homogeneous in the vector bundle structure. The terms $\pi_0$ and $c$ have degrees of homogeneity $-1$ and $-2$, respectively.
So,
\[
\Lie_X\pi = -\pi_0 - 2c .
\]
This Lie derivative is also a Poisson coboundary: $\Lie_X\pi = [X,\pi] = - \delta X$. So,
\[
-c = \pi - \delta X
\]
may be the cocycle we are looking for. 

If $p^*:H^2_{\mathrm{Lie}}(A) \to H^2_\pi(M)$ happens to be injective, then we can conclude that the quantization recipe for $(M,\pi)$ gives
\[
\cs_\P[\Sig(M),\sigma] \cong \cs[\Grp(A),\sigma_0]
\]
with any $\sigma_0$ such that $\Psi(\sigma_0)\simeq -c$. If $p^*$ is not injective, then this $\sigma_0$ is still correct, but the proof requires more analysis. I just give the construction here.

Consider the $\T$-extended groupoid $\G^{\sigma_0}$. Its Lie algebroid is an extension
\[
0 \to \R\times N \to \Alg(\G^{\sigma_0}) \to A \to 0
\]
classified by the cohomology class of $-c$. So, $\Alg(\G^{\sigma_0})$ can be identified with the bundle $A\oplus (\R\times N)$ with a bracket of the form \eqref{ext.bracket}. The dual $\Alg^*(\G^{\sigma_0})$ is a linear Poisson manifold which can be identified with $A^*\times \R$; the Poisson structure is essentially $\pi_0$ minus $c$ times the $\R$ coordinate. So, $M$ can be identified with the Poisson submanifold where the $\R$ coordinate equals $-1$.

The symplectic groupoid $T^* \G^{\sigma_0}$ integrates $\Alg^*(\G^{\sigma_0})$. The source and target of any element of this groupoid must lie in the same symplectic leaf, therefore they have the same $\R$ coordinate and this defines a map $H:\G^{\sigma_0}\to \R$.

By construction, there is a $\T$-action on $\G^{\sigma_0}$. This of course extends to a Hamiltonian $\T$-action on $T^*\G^{\sigma_0}$, where the Hamiltonian is just  $H$. The inverse image $H^{-1}(-1)\subset \G^{\sigma_0}$ is a subgroupoid over the submanifold $M$.
The symplectic reduction $H^{-1}(-1)/\T$ is the symplectic groupoid $\Sig(M)$. 
Let $Y\in\X^1[H^{-1}(-1)]$ be the restriction of the Hamiltonian vector field $\#dH$ generating this $\T$ action.

The Liouville form $\theta\in\Omega^1(T^*\G^{\sigma_0})$ is a symplectic potential, so the prequantization of $T^*\G^{\sigma_0}$ is given by the trivial bundle with connection $\nabla=d+i\theta$. From this, we can explicitly construct the prequantization $L\to\Sig(M)$. 
Because $Y\inner \theta = (\# dH)\inner \theta = H=-1$ on $H^{-1}(-1)$, the smooth sections of $L$ are identified with the functions $f\in\C^\infty(H^{-1}(-1))$ such that 
\[
0 = \nabla_Y f = Y(f) - if .
\]
 In other words, $L$ is the line bundle associated to the principal $\T$-bundle $H^{-1}(-1)\to \Sig(M)$ by the fundamental representation of $\T$.

Reducing this by the polarization gives the line bundle associated to $\G^{\sigma_0}\to\G$. In other words, $\sigma_0$ is the correct reduced twist.

\subsection{\Kahler-Poisson}
\label{Kaehler-Poisson}
\begin{definition}
A \emph{\Kahler-Poisson manifold} \cite{kar1} is a manifold equipped with both complex and Poisson structures, such that the Poisson bracket of any two local holomorphic functions vanishes.
\end{definition}

\begin{thm}
Let $\Sigma$ be a symplectic groupoid integrating a \Kahler-Poisson manifold $M$, and $\F\subset T_\co M$ the antiholomorphic tangent bundle. 
If 
\[
\P := T(\t,\s)^{-1}(\F\times\bar\F)
\]
is a regular distribution, then it is a polarization of $\Sigma$ as a symplectic groupoid.
\end{thm}
\begin{proof}
$\F$ gives an exact polarization $\F\times\bar\F$ of $\Pair(M)$. The groupoid anchor $(\t,\s):\Sigma\to M\times M=\Pair(M)$ is a homomorphism for any groupoid. By Lemma \ref{inverse.image}, $\P$ is a groupoid polarization.

The definition of \Kahler-Poisson implies that $\F^\perp$ is isotropic (with respect to $\pi$). Because $\t$ is a Poisson map, this implies that $T\t^*(\F^\perp)$ is isotropic; likewise, $T\s^*(\bar\F^\perp)$ is isotropic. These are subbundles of $T^{\s\perp}_\co\Sigma$ and $T^{\t\perp}_\co\Sigma$, which are symplectically orthogonal to one another. Therefore $\P^\perp=T\t^*(\F^\perp) + T\s^*(\bar\F^\perp)$ is isotropic, and $\P$ is coisotropic.

Again, use $\inv$ to split the exact sequence \eqref{A.sequence} and identify $\unit^*T\Sigma \cong TM\oplus \Alg(\Sigma) \cong TM\oplus T^*M$.
A straightforward computation shows that 
\[
\unit^*\P = \left\{(-\tfrac{i}2J\#\xi,\xi)\mid\xi\in T^*M\right\} 
\mbox,
\]
Where $J: TM\to TM$ is the complex structure.
So, $\rk \P = \dim M = \tfrac12\dim\Sigma$ and $\P$ is Lagrangian.
\end{proof}

Unlike the previous examples, this is not usually strongly admissible. Unless $\rk \pi$ is constant, $\P\cap\bar\P = T^\s_\co\Sigma\cap T^\t_\co\Sigma$ will not be a bundle.

It is rather difficult for this $\P$ to fail to be a regular distribution; indeed, I don't know if that ever happens. Such a failure at least requires a nonzero intersection $T\t^*\F^\perp \cap T\s^*\bar\F^\perp$, which is related to $\F\cap \ker \#$.

\begin{example}
Consider $S^2$ with any nonzero integrable Poisson structure.
A real polarization is impossible because $\P_0=0$ cannot be coisotropic, there exists no rank $1$ foliation of $S^2$, and $S^2$ is not parallelizable.
However, with any complex structure, this is \Kahler-Poisson. The intersection of $\F^\perp$ with $\ker\#$ is trivial, so this defines a polarization of any symplectic groupoid over $S^2$.
\end{example}

\subsection{Abelian and Locally Abelian}
In \cite{rie3}, Rieffel studied an explicit quantization for Poisson structures induced by Abelian group actions. In this section, I study this class of Poisson structures as well as the generalization where the Poisson structure is only locally of that form.  For these Poisson structures, I present a simple, explicit symplectic integration which is related to an action groupoid. I construct polarizations for these symplectic groupoids in two steps, using the results of Section~\ref{induced}.

Rieffel constructed an explicit deformed product. Because I have not studied quantization maps or deformed products here, I do not go so far as to compare my construction with his in any detail. 

Suppose that $T\subset\Diff(M)$ is a connected Abelian group of diffeomorphisms of a smooth manifold $M$. Let $\tf:=\Alg(T)$ be the Lie algebra of $T$.
Any translation invariant Poisson structure on $T$ is given by a bivector $\Pi\in\Wedge^2\tf$, and any such bivector determines an invariant Poisson structure. I use $\#_\Pi$ to denote the equivalent map, $\#_\Pi:\tf^*\to\tf$. 

Let $\action:\tf \to \X^1(M)$ be the action of $\tf$ by vector fields. If this is extended to exterior powers as $\action: \Wedge^\bullet\tf \to \X^\bullet(M)$, then
\beq
\label{Abelian_pi}
\pi := \action(\Pi)
\eeq
defines a $T$-invariant Poisson structure on $M$.
\begin{definition}
A Poisson structure induced from an Abelian group action in this way is an \emph{Abelian Poisson structure}.
\end{definition}

\subsubsection{Integration}
Another way to construct the Abelian Poisson structure \eqref{Abelian_pi} is to view $M$ as
\[
M \cong T \times_T M = (T\times M)/T .
\]
Then $\pi$ is just the push-forward of $\Pi\times 0$ by the quotient map. The algebra of continuous functions on $M$ is the $T$-invariant subalgebra 
\[
\C(M) \cong [\C(T)\otimes \C(M)]^T . 
\]
To quantize $(M,\pi)$ we just need to replace $\C(T)$ in this expression with an (equivariant) quantization of $(T,\Pi)$; see \cite{c-l}.

The symplectic integration of $T$ is $\Sig(T,\Pi) = \tf^*\ltimes T \cong T^*T$, where $\tf^*$ acts on $T$ by the composition of $\#_\Pi:\tf^*\to\tf$ and $\exp:\tf\to T$.  This is a group-groupoid, where the group structure is the Cartesian product of the Abelian group $\tf^*$ with $T$.

The quantization of $M$ should be (something like) the $T$-invariant subalgebra of the tensor product of $\C(M)$ with the quantization of $T$. The quantization of $T$ is constructed from the symplectic groupoid $\tf^*\ltimes T$, so the quantization of $M$ should be constructable from the groupoid
\[
(\tf^*\ltimes T)\times_T  M= \tf^*\ltimes_\Action M
\]
where the action $\Action$ is the composed map,
\[
\Action : \tf^* \xrightarrow{\#_\Pi} \tf \xrightarrow{\exp} T \subset \Diff(M) .
\]

In order to fit this into the framework of symplectic groupoids, we need to construct the symplectic integration of $M$ and find its relationship to $\tf^*\ltimes_\Action M$.

The dual of $\action:\tf\to\X^1(M)$ is a map, $\action^*:T^*M\to\tf^*$; equivalently, $\action$ extends to a Hamiltonian action on $T^*M$, and $\action^*$ is its momentum map. Let $p:T^*M\to M$ be the bundle projection. 
The anchor map for $T^*M$ factors through $\tf^*\times M$; for any $\xi\in T^*M$,
\[
\#\xi = \action_{\#_\Pi (\action^*\xi)}(p[\xi]) .
\]
With this in mind, define a map $q: T^*M\to \tf^*\times M$ by
\[
q(\xi)=(\action^*(\xi),p(\xi)) .
\]
This is a Lie algebroid homomorphism, if $\tf^*\times M$ is identified as the action Lie algebroid of $\tf^*$ on $M$, that is, the Lie algebroid of $\tf^*\ltimes_\Action M$.

Integration of Lie algebroids is functorial, so $q$ should integrate to a groupoid homomorphism from the symplectic integration $\Sig(M)$ to $\tf^*\ltimes_\Action M$. As it will turn out, $\Sig(M)$ is symplectomorphic to $T^*M$, so we can try to derive a groupoid structure on $T^*M$ by identifying $\tf^*\ltimes_\Action M$ with $\tf^* \times M$ and requiring that $q:T^*M\to\tf^*\ltimes_\Action M$ itself be a \emph{groupoid} homomorphism. 

Because $\tf^*$ is a vector space, we have the luxury of dividing by $2$, and can identify $\tf^*\ltimes_\Action M$ with $\tf^*\times M$ in a way that is more symmetrical than the standard presentation.
The source, target, multiplication, and inverse are:
\begin{gather*}
\s_{\tf^*\ltimes M}(u,x) = \Action_{- u/2}(x) \\
\t_{\tf^*\ltimes M}(u,x) = \Action_{u/2}(x) \\
(u,x)\cdot(v,y) = (u+v,\Action_{-v/2}x) = (u+v,\Action_{u/2}y) \\
(u,x)^{-1}=(-u,x)
\end{gather*}
where $u,v\in\tf^*$, $x,y\in M$.

Now, for $\xi\in T^*M$, the target must be $\t(\xi) = \t_{\tf^*\ltimes_\Action M}[q(\xi)] = \Action_{\action^*\xi/2}[p(\xi)]$. To simplify notation, define a map $\auxiliary : T^*M \to T\subset\Diff(M)$ by
\[
\auxiliary_{\xi} := \Action_{\action^*\xi/2} = \exp\left(\tfrac12 \#_\Pi\,\action^*\xi\right) .
\]
With this notation,
\[
\t(\xi) = \auxiliary_{\xi}[p(\xi)] \mbox,
\]
and likewise,
\[
\s(\xi) = \auxiliary_{-\xi}[p(\xi)] .
\]
The unit map of $\tf^*\ltimes_\Action M$ is the zero section, so the unit of $T^*M$ should be the zero section. The inverse of $(u,x)\in \tf^*\ltimes_\Action M$ is $(-u,x)$, so the inverse of $\xi\in T^*M$ should be $\xi^{-1}:=-\xi$.

Suppose that $\xi,\zeta\in T^*M$ are composable (i.e., $\s(\xi)=\t(\zeta)$) then 
\beq
\label{q_check}
q(\xi\cdot\zeta) = (\action^*\xi+\action^*\zeta,\s(\xi)) .
\eeq
 So, $\xi\cdot\zeta$ is a covector over $\s(\xi) = \auxiliary_{-\xi}[p(\xi)] = \auxiliary_{\zeta}[p(\zeta)]$ given by a sum of two terms. The natural candidate is to pull both $\xi$ and $\zeta$ back to this point and add them:
\[
\xi\cdot\zeta := \auxiliary^*_{+\xi}\zeta + \auxiliary^*_{-\zeta}\xi .
\]
 The map $\action^*$ is $T$-invariant, so this satisfies eq.~\eqref{q_check}.
\begin{thm}
With these source, target, unit, inverse, and product maps, the cotangent space $\Sigma := T^*M$ is a symplectic groupoid integrating $(M,\pi)$. The map $q:\Sigma\to\tf^*\ltimes_\Action M$ is a groupoid homomorphism.
\end{thm}
\begin{proof}
The groupoid axioms are straightforward to verify, and we have just verified the conditions for $q$ to be a homomorphism.

To check multiplicativity of the symplectic form, first observe that the manifold of composable pairs can be identified with the direct sum $T^*M\oplus T^*M$, by mapping a composable pair $(\xi,\zeta)\in\Sigma_2$ to the pair of cotangent vectors $(\auxiliary^*_{-\xi}\xi,\auxiliary^*_\zeta\zeta)$ over the point $\auxiliary_{-\xi}[p(\xi)]= \auxiliary_{\zeta}[p(\zeta)]$.

Let $\theta\in\Omega^1(T^*M)$ be the Liouville form. The symplectic form is multiplicative if and only if the simplicial coboundary $\partial^*\theta$ is closed. It is sufficient to check this over the dense submanifold where the action $\action$ has locally constant rank. We can use coordinates in which $\action(\tf)$ are constant vector fields and so $\pi$ is constant. A simple computation then shows that $\partial^*\theta$ is exact. It is the derivative of $-\frac12\pi$, viewed as a bilinear function on $T^*M\oplus T^*M$. So, $\partial^*\omega = -d[\partial^*\theta] =0$.

Finally, we need to verify that $\t:T^*M\to M$ is Poisson. It is again sufficient to check this over regular points of the action $\action$. 
\end{proof}


\subsubsection{Polarization}
As I have already mentioned, the symplectic integration 
\[
\Sig(T,\Pi) = \tf^*\ltimes T \cong T^*T
\] 
is a Lie group-groupoid (a special case of a double Lie groupoid). So, there are two ways of applying the Lie algebroid functor $\Alg$ to this. Applying $\Alg$ to the group structure gives a symplectic Lie algebra-groupoid, $\tf^*\oplus\tf$. An invariant polarization of $\tf^*\ltimes T$ is equivalent to a subspace $F\subset \tf^*_\co\oplus \tf_\co$ that is multiplicative, Hermitian, and Lagrangian.
Let 
\[
\daction : \tf^*\oplus\tf  \to \X^1(\tf^*\ltimes_\Action M)
\]
denote the product of the tautological action of $\tf^*$ on $\tf^*$ with the action $\action:\tf\to\X^1(M)$. This is a multiplicative action.
A polarization of $\tf^*\ltimes T$ extends trivially to a polarization of $(\tf^*\ltimes T)\times M$. This projects down to the polarization
\[
\Q := \Im \daction(F)
\]
of $\tf^*\ltimes_\Action M$, described in Lemma~\ref{push_forward}, provided that this is a regular distribution. Using Lemma~\ref{inverse.image}, the inverse image $Tq^{-1}(\Q)$ should at least be a groupoid polarization of $\Sigma$. 

\begin{prop}
If $F\subset \tf^*_\co\oplus\tf_\co$ is Lagrangian, Hermitian, and multiplicative, and $F$ is transverse to the isotropy subspace $\ker \action|_x \subseteq \tf$ for every $x\in M$, then 
\[
\P := Tq^{-1}(\Q) = Tq^{-1}[\Im \daction(F)]
\]
 is a symplectic groupoid polarization.
\end{prop}
\begin{proof}
First, note that $\ker \daction\vert_{(v,x)} = \ker \action\vert_x$. The transversality assumption means that $\daction$ maps $F$ injectively to each tangent fiber, therefore $\Q$ is regular and is a groupoid polarization of $\tf^*\ltimes_\Action M$.

Next we need to check that $\P\subset T_\co\Sigma$ is a regular distribution. It is easier to work in terms of the annihilator bundles where the definition of $\P$ is equivalent to,
\[
\P^\perp = Tq^*(\Q^\perp) \subset T^*\Sigma.
\]
 By $Tq^*(\Q^\perp)$, I mean the set of pulled-back covectors. The definition of $\Q$ can also be restated in terms of the annihilator;  $\Q^\perp$ is an inverse image, $\Q^\perp = (\daction^*)^{-1}(F^\perp)$.

The kernel of 
\[
Tq_\xi^*: T^*_{q(\xi)}(\tf^*\ltimes_\Action M) \to T^*_\xi\Sigma
\]
 is entirely vertical; it consists of the covectors $(v,0)$ where $v\in \tf^{**}=\tf$ such that $\action_v(p[\xi])=0$. The covector $(v,0)$ lies in $\Q^\perp$ if and only if $F^\perp \ni \daction^*(v,0) = (v,0)$ Since $F$ is Lagrangian, this is equivalent to  $(0,v)\in F$, and transversality then implies $v=0$.  So, $0 = \Q_{q(\xi)}\cap \ker Tq_\xi^*$ for every $\xi\in \Sigma$.  This shows that $Tq_\xi^*:\Q_{q(\xi)}^\perp\into T^*_\xi\Sigma$; Hence $\P^\perp$ --- and therefore $\P$ --- is regular.

The following diagram is commutative,
\[
\begin{CD}
T^*_\xi\Sigma @<{Tq_\xi^*}<< T^*_{q(\xi)}(\tf^*\ltimes_\Action M) @>{\daction^*}>> (\tf^*\oplus \tf)^* \\
@V{\#}VV @. @VV{\cong}V \\
T_\xi\Sigma @>{Tq_\xi}>> T_{q(\xi)}(\tf^*\ltimes_\Action M) @<{\daction\vert_\xi}<< \tf^*\oplus \tf \\
\end{CD}
\]
where the vertical map on the right is the symplectic identification of $\tf^*\oplus\tf$ with its dual;  either way around, $(v,\zeta)\in T^*_{q(\xi)}(\tf^*\ltimes_\Action M)$ maps to $(\action^*(\zeta),-\action_v[p(\xi)])$. The definitions of $\Q$ and $\P$ imply that this diagram restricts to,
\[
\begin{CD}
\P_\xi^\perp @<{Tq_\xi^*}<< \Q_{q(\xi)}^\perp @>{\daction^*}>> F^\perp \\
@V{\#}VV @. @VV{\cong}V \\
\P_\xi @>{Tq_\xi}>> \Q_{q(\xi)} @<{\daction\vert_\xi}<< F .\\
\end{CD}
\]
So, $\#\P^\perp\subseteq\P$, meaning that $\P$ is coisotropic.

The rank of $\Q$ is  $\dim F = \dim \tf$. The corank of $\P$ is the corank of $\Q$, which is 
\[
\dim (\tf^*\ltimes_\Action M) - \dim \tf = \dim M = \tfrac12\dim \Sigma \mbox,
\]
therefore $\P$ is Lagrangian.
\end{proof}

I derived this symplectic groupoid and polarization via the idea that the quantization would be given by $\tf^*\ltimes_\Action M$. 
Because the polarization $\P$ of $\Sigma$ is an inverse image polarization, the construction  of the algebra should factor through $q:\Sigma\to \tf^*\ltimes_\Action M$. The algebra ought to be a twisted polarized convolution algebra of $\Im q$, but $\Im q$ is not necessarily a Lie groupoid. It is not \emph{quite} all of $\tf^*\ltimes_\Action M$,
\[
\Im q = \{(v,x) \in \tf^*\ltimes_\Action M \mid v\in (\ker \action|_x)^\perp\}  .
\]
This has ``holes'' over points $x\in M$ where $\action|_x$ is not injective. 
The transversality assumption on $F$ implies that the polarization $\Q := \Im \daction(F)$ fills in these holes in $\Im q$. Any point of $\tf^*\ltimes_\Action M$ lies on a leaf of the singular foliation $\Q\cap\bar\Q$ which intersects $\Im q$.
So, the twisted, polarized convolution algebra of $\Im q$ should be that of $\tf^*\ltimes_\Action M$ after all.

\begin{example}
$\Im \#_\Pi\subset \tf$ is a symplectic vector space, so we can always choose a complex structure $J$ that makes this a \Kahler\ vector space. The subspace
\[
F := \{ (w,2iJ\#_\Pi w)\mid w\in \tf^*_\co\} \subset \tf^*_\co\oplus \tf_\co
\]
is Lagrangian, multiplicative, and Hermitian. This is transverse to $\tf_\co\subset \tf^*_\co\oplus \tf_\co$, so it always defines a polarization. So, the symplectic integration of an Abelian Poisson structure \emph{always} admits a polarization.
\end{example}

\begin{example}
If $F$ is real, then it is preserved by $\inv_{\tf^*\oplus\tf}$. The eigenspace decomposition shows that $F$ decomposes in the direct sum. Let $F_0:= F\cap \tf$ be the ``horizontal'' part. Since $F$ is Lagrangian, it must be $F = F_0^\perp \oplus F_0$. 

In this case, multiplicativity means simply that $\#_\Pi F_0^\perp \subseteq F_0$, i.e., $F_0$ is coisotropic. Transversality means that $F_0$ is transverse to each $\ker \action|_x$.
\end{example}

\subsubsection{Locally Abelian}
This class of examples can be generalized a little further. 
\begin{definition}
A Poisson structure on $M$ is \emph{locally Abelian} if it lifts to an Abelian Poisson structure on some covering manifold $\Mt$.
\end{definition}
This class of examples is important because I have shown in \cite{haw7} that if a compact Riemannian manifold admits a noncommutative deformation of its geometry, then the Poisson structure is locally Abelian with the Abelian group acting by isometries of $\Mt$.

Let $\Gamma$ be the covering group, so that $M=\Mt/\Gamma$. Obviously, the Poisson structure on $\Mt$ must be $\Gamma$-invariant, so $\Pi\in\Wedge^2\tf$ induces a well defined Poisson structure $\pi\in\X^2(M)$ if and only if it is $\Gamma$-invariant via the adjoint representation of $\Diff(\Mt)$. This implies that $T\subset \Diff(\Mt)$ can be chosen to be preserved by $\Gamma$.

The symplectic integration of $M$ can be obtained from the integration of $\widetilde M$ by $\Sig(\widetilde M) \cong T^*\widetilde M$. Symplectically, it is $\Sig(M) \cong T^*\widetilde M/\Gamma \cong T^*M$.

The (adjoint) action of $\Gamma$ extends to $\tf^*\oplus\tf$, and the construction of a polarization will work if and only if $F\subset\tf^*_\co\oplus\tf_\co$ is preserved by $\Gamma$. In particular, a real polarization is given by any $F_0\subset \tf$ that is coisotropic, transverse to the isotropy subspaces, and $\Gamma$-invariant.

\subsubsection{Quantization}
Consider a real polarization of a locally Abelian Poisson manifold, given by some $F_0\subset \tf$. For clarity, I only describe the case that $F_0$ is the Lie algebra of some (compact) torus $\T^k\subset T$.

In order for this polarization to be well defined, $F_0$ must be $\Gamma$-invariant. This means that $\T^k$ is a normal subgroup of the group of diffeomorphisms of $\widetilde M$ generated by $\T^k$ and $\Gamma$. So, the action of $\T^k$ on $\widetilde M$ descends to a well defined action of $\T^k$ on $M$ itself. 

By assumption, $F_0$ is transverse to the isotropy subspaces, this means that $\T^k$ acts freely on $\widetilde M$. There is some freedom in the choice of covering $\Mt$; we could replace $\widetilde M$ by $\widetilde M/(\T^k\cap \Gamma)$, and so we can assume that $\T^k\cap \Gamma$ is trivial. With this assumption, $\T^k$ acts freely on $M$, and so $M$ has the structure of a $\T^k$-bundle over the quotient manifold $N:=M/\T^k$. 
 
The quotient groupoid $\Sigma/\P$ is the quotient of $\tf^*\ltimes_\Action M$ by the action of $F=F_0^\perp\oplus F_0$; the two parts act independently on $\tf^*$ and $M$. The unit manifold is $M/F_0 = M/\T^k = N$. The first factor reduces to $\tf^*/F_0^\perp = F_0^* \cong \R^k$. The group action $\Action$ reduces to a well-defined action of $F_0^*$ on $N$. The quotient groupoid is thus $F_0^*\ltimes N$.

The leaves of this polarization have fundamental group $\Z^k$, so there are $k$ Bohr-Sommerfeld conditions. This is a generalization of the last example I gave in Section~\ref{trivial}, and by the same reasoning as there, the Bohr-Sommerfeld conditions reduce $F_0^*$ to $\Z^k$. So, the reduced groupoid is $\Z^k\ltimes N$.

Each point in $\Z^k$ corresponds to a character of $\T^k$, and again the twist line bundle over that component of $\Z^k\ltimes N$ is the line bundle constructed from the torus bundle $M\to N$ by that character.

\begin{example}
Suppose that $M$ is a trivial $\T^k$ bundle over $N$. Then the twist is trivial, and the quantization is $\cs(\Z^k\ltimes N) \cong \cs(\Z^k,\C_0(M))$. This is a case of the quantization constructed by Cadet \cite{cad1,cad2}. His general construction uses an open subgroupoid of $\Z^k\ltimes N$ in order to include cases when the action of $\T^k$ on $M$ is not free. 
\end{example}

\begin{example}
\emph{The Connes-Landi $3$-sphere.}
Consider a unit sphere $M=S^3$. A torus $T=\T^2$ acts on this by isometries, and any constant symplectic structure on $\T^2$ induces a Poisson structure on $S^3$. The two obvious generators of $\Alg(\T^2)$ do not act freely on $S^3$; their vector fields vanish on two great circles. We must choose $F_0$ transverse to these generators and the simplest choice is the diagonal.

So, we take $F_0\cong \R$ as the Lie algebra of $\T$ acting freely on $S^3$ with quotient $N = S^3 /\T \cong S^2$; $S^3$ is viewed as the tautological $\T$-bundle over $S^2$. The group $\Z$ acts on $S^2$ by rotations (the rotation angle is the symplectic area of $\T^2$). Lift this to some action of $\Z$ on the tautological line bundle. Let $U$ be the generator of this action on sections of the tautological line bundle.

The quantization of $S^3$ is the \cs-algebra generated by elements of the form $\psi U$ (and their adjoints) where $\psi$ is any smooth section of the tautological line bundle over $S^2$.

Note that changing the action on the line bundle would mean replacing $U$ with $uU$ for some function $u\in\C^\infty(S^2)$ with $\abs{u} =1$. This would give an isomorphic algebra.
\end{example}

\begin{example}
Heisenberg manifolds.
The $k$'th Heisenberg manifold is the total space of the $\T$-bundle over $\T^2$ with Chern class $\int_{\T^2}c_1 = k$. Rieffel \cite{rie2} constructed a quantization of this with two parameters $\nu,\mu\in \R$. 

The corresponding Poisson structure is locally Abelian. The symplectic foliation is the inverse image of the diagonal foliation of $\T^2$ with slope $\nu/\mu$; this is a Kronecker foliation if $\nu/\mu$ is irrational.

There is only one nice choice of polarization here. This is given by taking $F_0\cong\R$ to generate the principal action of $\T$ on $M$.

The reduced groupoid is $\Z\ltimes \T^2$, where the action of $\Z$ is generated by a $(\nu,\mu)$ translation. This can be lifted to an action on the line bundle associated to $M\to N$. Again, let $U$ be the generator of this action on sections.

The quantization of $M$ is the \cs-algebra generated by elements of the form $\psi U$, where $\psi$ is any smooth section of that line bundle.
\end{example}

\section{Outlook}
In this paper I have sketched a broadly ambitious programme for the geometric quantization of Poisson manifolds. Although I have only addressed the object side here, the central ambition is to realize quantization as a contravariant functor whose codomain is the category of \cs-algebras and $*$-homomorphisms.

The central idea is that the quantization of a Poisson manifold is a twisted polarized convolution \cs-algebra of a symplectic groupoid. Although I have only given a preliminary definition of this algebra,  I hope that the examples have convinced the reader that this idea unifies and extends the previously known geometric constructions of \cs-algebraic quantization. In particular, the cases where my preliminary definition applies do include the standard geometric quantization of symplectic manifolds, as well as twisted convolution algebras of $\s$-connected Lie groupoids. 

If quantization is a functor, then the domain category is not simply the category of Poisson manifolds. The main message of this paper is really about the objects of this category, i.e., the structure needed for quantization. I am proposing that an object in this category is a symplectic groupoid with a prequantization, polarization, and half-form structure.

By itself, this proposal raises many interesting questions, because any construction or question in Poisson geometry can be revisited by incorporating some or all of this additional structure.

Not all Poisson manifolds should be quantizable, and the quantization of a Poisson manifold is not generally unique.
In each step of this recipe for quantization, there are existence and uniqueness issues.

First, the Poisson manifold $M$ needs to be integrable. The integrability conditions were found by Crainic and Fernandes \cite{c-f2}. 
In general, there exists more than one symplectic groupoid $\Sigma$ integrating $M$; however, any choice of $\Sigma$ is a quotient of  the unique $\s$-simply connected integration $\Sig(M)$. This nonuniqueness seems to account for the difference between standard geometric quantization of a symplectic manifold, and my construction in \cite{haw5}.

Second, a prequantization of the symplectic groupoid $\Sigma$ may not exist. Theorem~\ref{c_prequantizability} (due to Crainic \cite{cra2}) describes the prequantizability condition for $\Sig(M)$ and shows that its prequantization is unique. In general, $\Sigma$ is prequantizable if $\Sig(M)$ is prequantizable and the prequantization descends to $\Sigma$. Prequantization of $\Sigma$ is not generally unique, but this nonuniqueness is described by $H^1(\Sigma,M;\T)$ according to Theorem~\ref{nonunique}.

Third, a polarization may not exist and may not be unique. Understanding these issues requires a general description of the polarizations of Poisson manifolds.

Fourth, the construction of a half-form bundle (or sheaf) $\Omega^{1/2}_\P$ involves taking the square root of a line bundle. This may not exist and may not be unique.

Finally, there may be some freedom in choosing the completion to a \cs-algebra.

This paper is only intended as a beginning. There are several fundamental questions that remain to be answered in order to carry this programme of quantization forward.

What is a polarization of a Lie algebroid or Poisson manifold? I presented here a description of \emph{real} groupoid (and symplectic groupoid) polarizations at the infinitesimal level, in terms of the Lie algebroid. This needs to be extended to complex polarizations. 

What is a twisted polarized \cs-algebra for a symplectic groupoid? I have given a preliminary definition of such a convolution algebra in the best behaved cases. I have not discussed the details of completing this to a \cs-algebra, but as I suggested in Section~\ref{Completion}, there should be \emph{at least} two ways of doing this. I suggested that the general twisted, polarized convolution algebras may be constructed from a total sheaf cohomology space. This leaves the problem of defining convolution on cohomology. This may be extremely difficult, because it is closely related to the problem of defining an inner product on cohomological wave functions in conventional geometric quantization, which has never been resolved.

To what extent is this algebra independent of the polarization? In the examples where I have considered more than one polarization of the same groupoid, the resulting algebras were isomorphic. This happened for a symplectic vector space, a torus, and manifolds with zero Poisson structure. I doubt that the algebra will always be independent of the polarization, but it would be useful to characterize which polarizations are equivalent in this way. 

How are Poisson maps quantized? A smooth map intertwining Poisson brackets is the semiclassical analogue of a homomorphism of noncommutative algebras. Such a Poisson morphism should be quantizable to a $*$-homomorphism when it is sufficiently compatible with the polarizations. A construction for quantizing some Poisson maps would make quantization a functor.

Is this really quantization? It is important (or at least desirable) to prove that this construction satisfies some reasonable definition of \cs-algebraic deformation quantization along the lines sketched in Section~\ref{Quantization}.

\section*{Acknowledgements}
In the course of this work so far, I have benefited from discussions with Klaas Landsman, Marius Crainic, Hessel Posthuma, Ieke Moerdijk, and Alan Weinstein. This work was supported by a Marie Curie Incoming International Fellowship from the European Commission.

\end{document}